\documentclass[journal]{IEEEtran}   
\pdfoutput=1
\usepackage{cite}
\usepackage{amsmath,amsthm}
\usepackage{array}
\usepackage{amssymb}
\usepackage{bm,bbm}
\usepackage{fixltx2e}
\usepackage{enumerate}     

\usepackage{microtype}     
\usepackage{hyperref}      
\usepackage{longtable}

\usepackage{ifpdf}

\ifCLASSINFOpdf
   \usepackage[pdftex]{graphicx}
   \graphicspath{{../pdf/}{../jpeg/}}
\else
   \usepackage[dvips]{graphicx}
\fi

\usepackage{subfig}
\usepackage{url}

\usepackage{setspace}

\usepackage{cases}

\newcommand{\km}{{k_{c}}}

\newcommand{\ra}{\rightarrow}

\newcommand{\rmx}{{\rm x}}
\newcommand{\e}{{\eta_{1}}}
\newcommand{\ec}{{\eta_{1,c}}}

\newcommand{\om}{\omega}

\newcommand{\bfs}{{\bf s}}
\newcommand{\bfk}{{\bf k}}
\newcommand{\kk}{k}
\newcommand{\bfr}{{\bf r}}

\newcommand{\rr}{r}

\newcommand{\EE}{\mathbb{E}}
\newcommand{\Rr}{\mathbb{R}}

\newcommand{\inc}{\Gxx}

\newcommand{\Gxx}{C_{\rmx\rmx}}

\newcommand{\D}{\mathcal{D}}

\newcommand{\Fr}{\mathcal{F}}
\newcommand{\iFr}{\mathcal{F}^{-1}}
\newcommand{\A}{\mathcal{A}}

\newcommand{\bmthe}{{\bm \theta}}

\theoremstyle{plain}
\newtheorem{theorem}{{\bf Theorem}}[section]

\newtheorem{definition}{{\bf Definition}}[section]
\newtheorem{proposition}{Proposition}

\begin{document}

\title{Radial Covariance Functions Motivated by Spatial Random Field Models with Local Interactions}


\author{Dionissios~T.~Hristopulos\textsuperscript{*}%
    \thanks{D.~T.~Hristopulos is with the School of Mineral Resources Engineering, Technical University of Crete, Chania 73100, Greece, (e-mail: dionisi@mred.tuc.gr).}
    }
\markboth{}{Hristopulos: Covariance Models based on Local Interactions}

\maketitle

\begin{abstract}
Random fields based on energy functionals with local interactions
possess flexible covariance functions, lead to computationally efficient algorithms for spatial data processing, 
and have important applications in Bayesian field theory.
We derive explicit expressions for a family of
 radially symmetric, non-differentiable,  Spartan covariance functions in $\Rr^2$ that involve
 the modified Bessel function of the second kind.
In addition to the characteristic length and the amplitude coefficient,  the Spartan
covariance parameters include
the rigidity coefficient $\e$ which determines the shape of the covariance function.
If $ \e \gg 1$ Spartan covariance functions exhibit multiscaling.
We also derive a  family of
radially symmetric, infinitely differentiable Bessel-Lommel covariance
 functions valid  in $\Rr^{d}, d\ge 2$.
 We investigate the  parametric dependence of the integral range
 for Spartan and Bessel-Lommel covariance functions using explicit relations and numerical simulations.
 Finally, we define a generalized spectrum of correlation scales $\lambda^{(\alpha)}_{c}$ in terms of the fractional Laplacian of the covariance function;
 for $0 \le \alpha \le1$ the  $\lambda^{(\alpha)}_{c}$
   extend from the smoothness microscale $(\alpha=1)$ to the  integral range $(\alpha=0)$.
   The smoothness scale of mean-square continuous but non-differentiable random fields vanishes;
    such fields, however, can be discriminated by means of  $\lambda^{(\alpha)}_{c}$ scales obtained for
    $\alpha <1$.
\end{abstract}



\section{Introduction}
\subsection{Background}
\IEEEPARstart{T}{he} theory of  \emph{spatial random fields (SRFs)} is
a powerful mathematical framework for modelling spatial variability~\cite{Adler81,Yaglom87,Vanmarcke10}.
The SRF theory has a multidisciplinary scope of
applications that include fluid mechanics~\cite{Monin71}, computational and probabilistic engineering mechanics~\cite{Spanos03},
materials science~\cite{Ostoja98,Torquato02,Sena13}, hydrological modeling~\cite{Gelhar93,Kitanidis97,Rubin03},
petroleum engineering~\cite{Hohn99,Deutsch02,Farmer05}, environmental monitoring
\cite{Christakos92,Smith00,Kanevski04}, mining exploration and mineral reserves
estimation~\cite{Goovaerts97,Armstrong98}, environmental health~\cite{ch98},
geophysical signal processing~\cite{Smith05},
image analysis~\cite{Winkler95,Wilson02}, machine learning~\cite{Rasmussen06}~\footnote{The term \emph{Gaussian processes} is used for SRFs in the machine learning literature.} statistical cosmology~\cite{Bertschinger01,Peebles93},
medical image registration
\cite{Alzola05}, as well as in structural and functional mapping of the brain~\cite{Siegmund95,Leow04,Cao01}.

A second-order description of Gaussian or non-Gaussian, centered, scalar SRFs  requires the introduction of
a permissible covariance function, i.e., a non-negative definite, two-point function~\cite{Bochner59,Adler81}.
Covariance functions are the key ingredient in Best Linear Unbiased Prediction (BLUP) methods for
 spatial interpolation and in geostatistical simulation~\cite{Cressie93,Stein99,Lantuejoul02,Wackernagel03}.
Covariance functions are also used as non-negative definite kernels
in machine learning approaches based Gaussian processes~\cite{Rasmussen06,Hofmann08} as well as in radial basis function interpolation~\cite{Ohtake04,Wendland05}
and surface reconstruction~\cite{Dinh02}.
Valid covariance functions are
employed in the simulation of non-Gaussian random fields via the
Karhunen-Lo{\`e}ve expansion~\cite{Sakamoto02,Phoon05}.
Most covariance  functions used in the analysis of spatial data are mathematical constructions
that satisfy Bochner's theorem~\cite{Bochner59}. In principle, covariance functions can also be obtained from the solutions of {stochastic partial
differential equations}~\cite{ch98,Kolovos04,Lim09} or from the solution of statistical field theories~\cite{Kardar07}.
It is often impossible, however,
to obtain closed-form covariance solutions for the entire distance range based on these approaches.
The ``standard'' covariance models (e.g., the exponential, spherical,  and Gaussian functions) involve  two parameters: the variance and the
correlation length. The variance determines the amplitude of the fluctuations, whereas all
 length scales characterizing SRF patterns are trivially proportional to the correlation length.
Hence, two-parameter covariance models may fail to capture the variability of spatially distributed processes.
Whereas some models with more than two parameters exist in the literature, e.g.
 the rational quadratic~\cite{Genton02,Rasmussen06}, the Whittle-Mat\'{e}rn~\cite{Whittle54,Stein99},
 additional flexible covariance models are needed~\cite{Genton02}.

\subsection{Approach and Main Contributions}

Our approach is based on SRFs governed by  a
 Gibbs energy functional with local self-interactions~\cite{dth03b}.
 The so-called \emph{Spartan spatial random fields (SSRFs)}
 have a rational spectral density that involves a fourth degree polynomial in the denominator, c.f. equation~\eqref{eq:ssrf-spd} below.
The energy functional is not meant to represent an actual energy of SRF configurations but to provide
\emph{ad hoc} geometric constraints that render SRF states with certain spatial patterns more likely than others.
In one and three dimensions, explicit SSRF covariance expressions
  were derived at the limit of infinite
 \emph{ultraviolet (high-wavenumber)} cutoff of the SSRF spectral density function~\cite{dthsel07}.

 The main contributions of this paper are threefold. First,
   we derive explicit expressions for stationary, non-differentiable, three-parameter, radial SSRF covariance functions
  $\bfr \in \Rr^2 \mapsto \Gxx(\bfr) \in \Rr$, and we show that
  these functions exhibit multiscaling in a specified region of the parameter space.

Secondly, we obtain  stationary, infinitely differentiable,  radial Bessel-Lommel covariance functions
 for $\bfr \in \Rr^d, \, d \ge 2$ based on the
 reciprocal SSRF spectral density  with finite  spectral cutoff.
 These functions  generalize the oscillatory ordinary Bessel covariance functions and
 provide covariance models for smooth spatial processes in $\Rr^2$ and $\Rr^3$.
They can also be used  in machine learning applications that involve higher-dimensional
 feature spaces~\cite{Genton02,Rasmussen06}.

 Thirdly, we define a generalized correlation spectrum that captures the integral ranges of
 fractional SRF derivatives.
 The spectrum can be used with both differentiable and non-differentiable stationary SRFs.
 In contrast with the smoothness microscale, which vanishes for mean square continuous but
 non-differentiable SRFs, the proposed spectrum discriminates between
 non-differentiable SRFs with different covariance functions.

\subsection{Structure}
 The remainder of this paper is structured as follows: Section~\ref{sec:prelim} presents
 notation and mathematical background.
Section~\ref{sec:spartan} is a brief, non-mathematical overview of SSRFs. Section~\ref{sec:SSRF-2d}
derives two-dimensional SSRF
covariance functions at the limit of infinite ultraviolet cutoff of the spectral density.
In Section~\ref{sec:BL}, Bessel-Lommel covariance functions
 in  $d\ge 2$ are obtained. Section~\ref{sec:cor-scales} defines a generalized correlation spectrum and investigates
 its  dependence for
 SSRF and Bessel-Lommel SRFs.
 Section~\ref{sec:simul} investigates geometric features
 of Gaussian SRF realizations  with SSRF
 and Bessel-Lommel covariance functions in connection  with the results of Section~\ref{sec:cor-scales}.
 Section~\ref{sec:concl} summarizes the results and discusses potential applications.
Most mathematical proofs are relegated to Appendices.
Finally, Appendix~\ref{app:table} contains a table that summarizes the notation used.

\section{Definitions and Preliminaries}
\label{sec:prelim}

Given a complex number $z$, $\bar{z}$  denotes its complex conjugate;
$\Re(z)$ and $\Im(z)$ denote respectively the real and imaginary parts
of $z$, i.e., $z=\Re(z) + \jmath \, \Im(z)$. An operator $\mathbb{A}$ is
self-adjoint if $\mathbb{A}=\mathbb{A}^{\dag}$, where $\mathbb{A}^{\dag}$ is the
complex conjugate of $\mathbb{A}.$ The transpose of a vector or matrix ${\bf V}$ will be denoted by
${\bf V}^T$.

A scalar SRF
 $ \{ X(\bfs,\om) \in \mathbb{R}; \, {\bf s} \in \mathcal{D}; \omega \in \Omega \}$ is a mapping
from the probability space $(\Omega,A,P)$ into the space of real
numbers so that for each fixed ${\bf s}$,
$ X(\bfs,\om)$ is a measurable function of $\omega$ \cite{Adler81}.
The spatial domain $\mathcal{D}$   includes all points ${\bf s} \in  \mathbb{R}^{d}$.
The reciprocal space involves the wavevectors $\bfk \in \Rr^{d}.$
In addition to $d$, we use as needed the \emph{dimension index} $\nu=d/2-1$ to simplify notation.
An SRF involves by definition a multitude of probable \emph{states} or
\emph{realizations} indexed by $\omega$~\cite{Yaglom87,Christakos92}.
If $\Phi: \Rr \ra \Rr$ is a Borel measurable function, then
$\EE\left[\Phi\left(X(\bfs,\om)\right)\right]$ denotes the \emph{expectation} of $\Phi\left(X(\bfs,\om)\right)$
 over the ensemble of states $\Omega$.
A \emph{sample (realization)} of the SRF will be denoted by $x(\bfs) \in \mathbb{R}.$

An SRF is \emph{second-order stationary} or \emph{weakly stationary}
if (i) its expectation $\EE[X(\bfs,\om)]=m_{\rmx}$ is independent of the
location and (ii) its autocovariance function
depends only on the spatial lag vector ${\bf r} \in \Rr^{d}$~\cite[pp.
308-438]{Yaglom87}, i.e.,
\begin{align}
\label{eq:cov-def}
\Gxx(\bfr) &  = \EE[X(\bfs,\om)  X(\bfs+\bfr,\om)]
  - m_{\rmx}^{2}.
\end{align}
In the following, we use the term ``stationary'' to refer to \emph{second-order
stationarity}.
A stationary, scalar SRF is \emph{statistically isotropic} (for short isotropic) if its covariance function
is \emph{a radial function}, i.e., if
$\Gxx(\bfr)=\Gxx(\|\bfr\|)$, where $\|\bfr\|$ is the Euclidean norm of
 $\bfr$. In the following, we will use $\rr$ and $\kk$ to denote the Euclidean norms
 of the vectors $\bfr$ and $\bfk$ respectively~\footnote{This notation is ambiguous in $d=1$ but
 well defined for $d \ge 2$.}.

\subsection{Spectral Representation}
\label{ssec:ssrf-spectral}
We review the spectral representation of  covariance functions, since we will use it below
to derive covariance functions in position space.
For a stationary SRF let the \emph{Fourier transform} of the covariance function from the space of lag vectors
$\bfr \in \Rr^d$ to the
reciprocal (Fourier) space of wavevectors $\bfk \in \Rr^d$ be defined by
means of the following improper multiple integral:
\begin{equation}
    \label{eq:covft}
    \widetilde{\Gxx}({\bf k};\bmthe) = \Fr[ \Gxx({\bf r};\bmthe)] = \int_{\Rr^d} d{\bf r}\; {\rm e}^{-\jmath \bf{k\cdot r}}
    \Gxx({\bf r};\bmthe),
\end{equation}
where  $d{\bf r} = \prod_{i=1}^{d} dr_{i}$.   The \emph{inverse Fourier transform}
is then given by means of the improper multiple integral
\begin{equation}
    \label{eq:invcovft}
    \Gxx({\bf r};\bmthe) = \iFr[ \widetilde{\Gxx}({\bfk};\bmthe)] = \frac{1}{(2\,\pi)^{d}}\,\int_{\Rr^d} d{\bf k}\;
    {\rm e}^{\jmath \bf{k\cdot r}} \, \widetilde{\Gxx}({\bf k};\bmthe).
\end{equation}
The integrals~\eqref{eq:covft} and~\eqref{eq:invcovft} are well defined if $\Gxx({\bf r};\bmthe)$ and $\widetilde{\Gxx}({\bf k};\bmthe)$
are absolutely integrable, respectively.

\begin{theorem}{Bochner-Khinchin's permissibility theorem}
\label{theor:Bochner}
 requires that
\begin{enumerate}
\item $\widetilde{\Gxx}(\kk;\bmthe) >0, \, \forall \bfk$,

\item
  $\int_{\Rr^d} d\bfk \,  \widetilde{\Gxx}(\kk;\bmthe) < \infty$,

 \end{enumerate}
 for the inverse Fourier transform~\eqref{eq:invcovft} to
represent a valid covariance function~\cite{Bochner59}, \cite[p.106]{Yaglom87}.
\end{theorem}

For \emph{radial covariance} functions the pair of
 transforms~\eqref{eq:covft} and~\eqref{eq:invcovft}  is expressed as follows~\cite[p. 353]{Yaglom87}
\begin{subequations}
\begin{equation}
    \label{eq:covft-iso}
    \widetilde{\Gxx}({\kk};\bmthe)= \frac{(2\pi)^{d/2}}{\kk^{\nu}}  \int_{0}^{\infty} d{\rr} \, \rr^{d/2}
    {J_{\nu}(\kk  \rr)}  \Gxx({\rr};\bmthe),
\end{equation}
\begin{equation}
    \label{eq:invcovft-iso}
    \Gxx({\rr};\bmthe)= \frac{1}{(2\pi)^{d/2} \rr^{\nu}} \int_{0}^{\infty} d{\kk} \, \kk^{d/2}
    {J_{\nu}(\kk  \rr)}     \widetilde{\Gxx}({\kk};\bmthe),
\end{equation}
\end{subequations}
where $\nu = d/2-1$, $J_{\nu}(x), \, x\in \mathbb{R},$ is the \emph{Bessel function of the first kind of
order} $d/2-1,$ and $\kk$ is the Euclidean norm of the wavevector $\bfk$~\cite{Schoenberg38}.

\section{Review of Spartan Spatial Random Fields}
\label{sec:spartan}
Gibbs random fields have a joint probability density function defined in terms of an energy functional
of the sample function $x(\bfs)$.
Spartan Spatial Random Fields (SSRFs) are Gibbs random fields whose energy functional is
defined in terms of local interactions between the values of
$x(\bfs)$ at neighboring points.
A particular type of SSRF is the Fluctuation-Gradient-Curvature FGC-SSRF~\footnote{Henceforward SSRF for simplicity.} model
whose energy functional involves the square gradient and the square Laplacian of $x(\bfs)$~\cite{dth03,dthsel07}.

The SSRF model is determined by the parameter vector
$\bmthe=(\eta_0, \e, \xi, \km)^T,$ which includes the \emph{amplitude coefficient} $\eta_0$,
the \emph{characteristic length} $\xi$,  the \emph{rigidity coefficient} $\e$,
and  the \emph{spectral cutoff}, $\km$.
The amplitude coefficient is analogous to temperature in statistical mechanics and controls
  the overall magnitude of the fluctuations.
  The rigidity coefficient controls the resistance of the field to changes of its gradient.
  The characteristic length controls, in conjunction with $\e$, the relative strength of the   square Laplacian versus the square gradient term.
  The spectral cutoff represents an implicit upper bound
 in reciprocal space, which should be finite for the SSRF states to be mean square differentiable (see Proposition~\ref{prop:ssrf} and its proof).
 We derive explicit covariance functions at the limit $\km \ra \infty$ where the
 SSRF states become non-differentiable in the mean square sense. These functions represent
 asymptotic limits of the covariance functions that
 correspond to the SSRF energy functional; they are nevertheless permissible, albeit non-differentiable, covariance functions.

 The \emph{SSRF spectral density} is given by the following function as shown in~\cite{dth03b}
\begin{subequations}
\label{eq:ssrf-spd}
\begin{equation}
 \label{eq:covspd-1}
    \widetilde{\Gxx}({\bf k};\bmthe)=
    \frac{\eta _0 \,\xi ^d \, \mathbbm{1}_{\km \ge \kk}(\kk)}{\Pi(\kk \xi)},
\end{equation}
\begin{equation}
\label{eq:char-pol}
\Pi(u) =1 + \e u^2 + u^4,
\end{equation}
\end{subequations}
where  $\Pi(u)$, $u=\kk \, \xi$,
is the \emph{characteristic polynomial} of the SSRF spectral density,  and $\mathbbm{1}_{A}(x)$ is the
\emph{indicator function}
of the set $A$, i.e., $\mathbbm{1}_{A}(x)=1, \, x \in A$
and $\mathbbm{1}_{A}(x)=0, \, x \ni A.$

The non-negativity of $\widetilde{\Gxx}({\bf k};\bmthe)$ ---Condition (1) of Bochner's theorem--- is
 ensured $\forall \km >0$, if $\e> -2$~\cite{dth03b}.
 The integrability, Condition (2)  of Bochner's theorem,
 is satisfied for  $\forall d\ge 1$, $ \forall \km \in \Rr_{+}$, or if $\km \ra \infty$  for $d<4$.

For radial SSRF covariance functions, the \emph{spectral representation}   is given by
 the following one-dimensional integral derived from~\eqref{eq:invcovft-iso}
\begin{equation}
\label{eq:cov-spectral} \Gxx({\rr };\bmthe) = \frac{\eta _0\,\xi
^d \,\rr^{1-d/2} }{(2\pi)^{d/2} } \int_{0}^{\km} d \kk
 \frac{\kk^{d / 2} J_{d/2-1}(\kk \rr)}  {{1 + \eta _1 (\kk\xi )^2  + (\kk\xi )^4 }}.
\end{equation}

The spectral integral~\eqref{eq:cov-spectral} yields  functions $\Gxx({\bf r};\bmthe) $ in position
space that change with $d$.
In $\Rr^1$ and $\Rr^3$  the evaluation of the spectral integral~\eqref{eq:cov-spectral}  is facilitated
by the fact that $J_{-1/2}(\cdot)$ and $J_{1/2}(\cdot)$
can be expressed in terms of trigonometric and
algebraic functions respectively. The integral is then evaluated at the limit $\km \ra \infty$ using Cauchy's Theorem of Residues~\cite{dthsel07}.
Nonetheless, many ``signals'' of interest, including digital images and
data from environmental sensor networks are strictly or
approximately two-dimensional fields.

\section{SSRF Covariance in Two Dimensions}
\label{sec:SSRF-2d}

Below we derive explicit
expressions for SSRF covariance functions
in $\Rr^2$ at the limit $\km \ra \infty$.
The spectral integral~(\ref{eq:cov-spectral})
is expressed by the following \emph{Hankel transform}
of order zero:
\begin{equation}\label{eq:cov2d-spectral}
  \Gxx({\bf r};\bmthe ) = \frac{\eta_0 \xi^2 }{2\pi} \,
  \int\limits_{0}^{\infty } d\kk\, \frac{ \kk \, J_{0}( \kk \rr)}{1+\e (\kk\xi)^2+ (\kk \xi)^4}.
\end{equation}
Hankel transforms of integer order often do not admit explicit expressions.
Nevertheless, the integral~\eqref{eq:cov2d-spectral} is amenable to explicit evaluation by means of
 the Hankel-Nicholson integration formula~\cite[p.~488]{Abramowitz72}.

\begin{proposition}
\label{prop:ssrf}
The  SSRF covariance function in $\Rr^2$ defined by the integral~\eqref{eq:cov2d-spectral} for $\e > -2 $
is given by the following equations
\begin{subequations}
\label{eq:cov2d-ssrf}
\begin{equation}
\label{eq:Cov_eta_gt_2}
\Gxx(h ;\bmthe) =\frac{\eta_0 \left[  K_{0} (hz_{+})   -  K_{0} (hz_{-}) \right] }{2\pi \sqrt{\eta_{1}^{2}-4}}
, \; \e > 2
\end{equation}
\begin{equation}
\label{eq:Cov_eta_eq_2}
       \Gxx(h;\bmthe) = \left( \frac{\eta_0 \, h}{4\pi} \right) \, K_{-1}(h), \;  \e = 2.
\end{equation}
\begin{equation}
\label{eq:Cov_eta_lt_2}
\Gxx(h;\bmthe)  = \frac{\eta_0  \, \Im\left[  K_{0} (h z_{+}) \right]}{\pi \sqrt{4 - \eta_{1}^{2}}},  \;
 |\e| <2,
\end{equation}
\label{eq:cov2d}
\end{subequations}
where $h=\rr/\xi$  is the dimensionless two-point lag distance,
 $K_{\nu}(\cdot)$ is the modified Bessel function of the second
kind and order $\nu$,  $\Im\left[\cdot\right]$ denotes the imaginary part, and
\begin{equation}
\label{eq:zpm}
z_{\pm}^2={\left(\e \mp \Delta\right)}/{2}.
\end{equation}
SSRFs  with  covariance functions given by~\eqref{eq:cov2d-ssrf} are continuous but non-differentiable in the mean square sense.
\end{proposition}

\begin{IEEEproof}
The proof is given in the Appendix~\ref{App:proof-ssrf}.
\end{IEEEproof}

For $\e =2$, Eq.~\eqref{eq:Cov_eta_eq_2}  recovers the known \emph{K-Bessel covariance function}~\cite{Lantuejoul02}.

\begin{figure}[!t]
  \centering
    \includegraphics[width=0.85\linewidth]{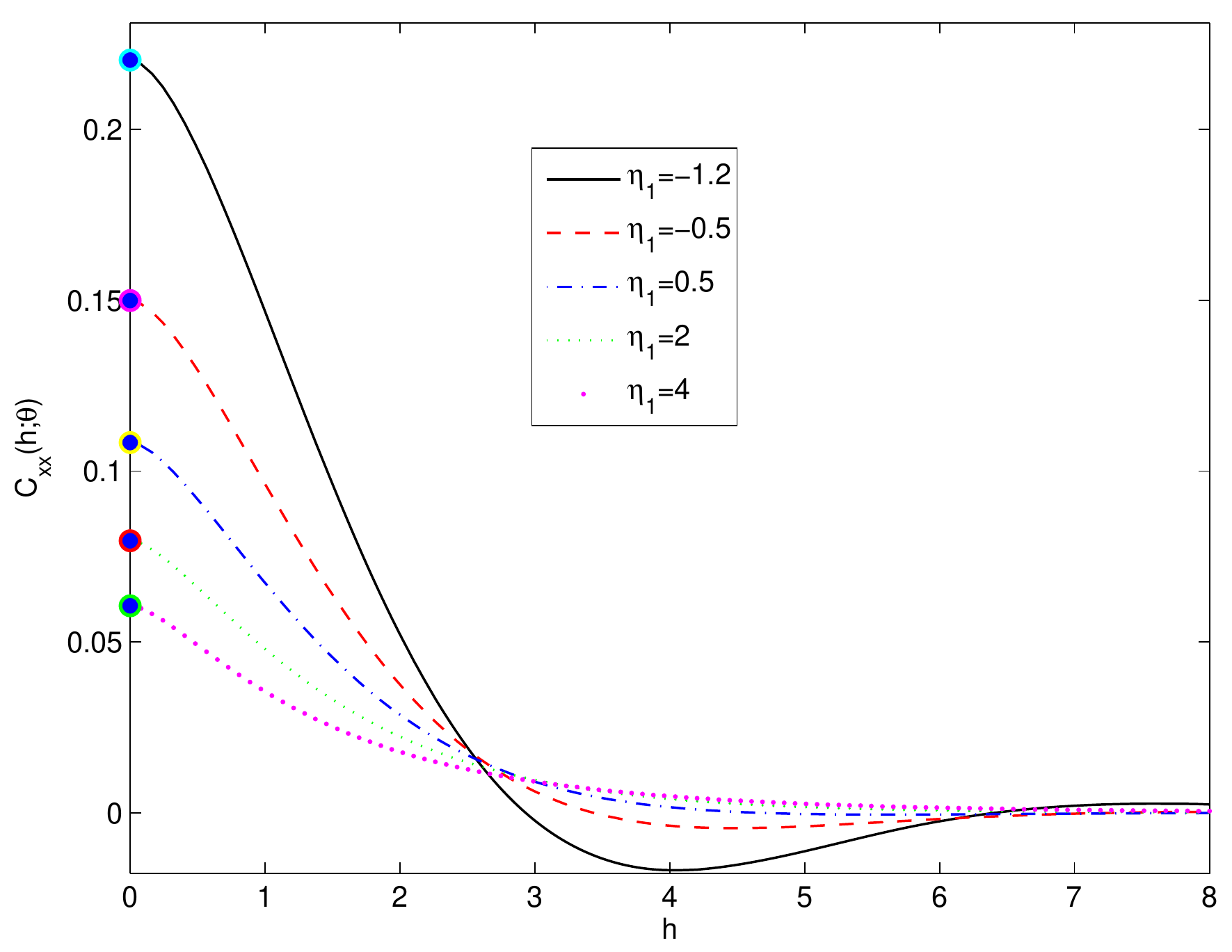}
  \caption{Covariance function $\Gxx(h;\bmthe)$ versus normalized distance $h$ for different values of $\e$.
  For $\e>2$ $\Gxx(h;\bmthe)$ is calculated from Eq.~\eqref{eq:Cov_eta_gt_2},
  for $\e=2$, $\Gxx(h;\bmthe)$ is calculated from Eq.~\eqref{eq:Cov_eta_eq_2},
   whereas for $|\e|<2$, $\Gxx(h;\bmthe)$ is calculated from Eq.~\eqref{eq:Cov_eta_lt_2}.
   In all cases, $\eta_0=1$. Filled circles at the origin $(h=0)$ represent the independently calculated variance.
   }
  \label{fig:cov_2d}
\end{figure}

We illustrate the dependence of $\Gxx(h;\bmthe)$ on $h$ for different  $\e$
in Fig.~\ref{fig:cov_2d}.  For $\e <0$,
$\Gxx(h;\bmthe)$ exhibits a negative peak and then returns to positive values; this damped
oscillatory behavior is due to the peak of the spectral density  for $-2<\e<0$~\cite{dth03b}.
The negative peak becomes more pronounced as $\e \ra -2$.
The values  $\Gxx(h=0;\bmthe)$
 agree with the  variance expressions
which were independently obtained in~\cite{dth03b}.
In general, higher values of $\e$ reduce the variance and increase the spatial coherence
 as evidenced in the slower decay of the tail of $\Gxx(h;\bmthe)$. This behavior reflects the higher stiffness of SSRF realizations as $\e \uparrow$ and is
 further elaborated in Section~\ref{sec:cor-scales}.

\section{Bessel-Lommel Covariance Functions}
\label{sec:BL}
We construct covariance functions for SRFs with spectral density proportional
to $\Pi(u)$ given by~\eqref{eq:char-pol} within a finite
spectral band which is cut off at $\km$.
The explicit expressions  involve products of  Bessel functions of the first kind and
 Lommel functions. This new class of
 non-negative definite, infinitely differentiable  functions  is  suitable for  smooth
 spatial processes defined  in Euclidean  $\Rr^2, \Rr^3$ spaces. In addition, they provide
   flexible kernel functions for
  machine learning applications in high-dimensional spaces~\cite{Genton02}.

Based on the general expression~\eqref{eq:invcovft-iso}
for the radial inverse Fourier transform, the Bessel-Lommel covariance function
$\inc^{BL}(\bfr;\bmthe)$ is  given by the following spectral integral
\begin{align}
\label{eq:cov-BL-spectral}
\inc^{BL}(\bfr;\bmthe) & = \frac{\rr}{(2\pi  \rr)^{d/2} } \int_{0}^{\km} d \kk \,
 {\kk^{d / 2} J_{d/2-1}(\kk \rr)} \,  \widetilde{\inc^{BL}}(\kk;\bmthe).
\end{align}
The spectral density  $\widetilde{\inc^{BL}}(\kk;\bmthe)$  is given by
\begin{equation}
\label{eq:spd-BL}
\widetilde{\inc^{BL}}({\kk ; \bmthe})= c_0 + c_{1} \kk^2 + c_{2} \kk^4,
\end{equation}where $c_0 = (\eta_0 \xi^d)^{-1},$ $c_{1} = \e \,(\eta_0 \xi^d)^{-1} \xi^2,$ and
$c_{2} = (\eta_0 \xi^d)^{-1} \xi^4.$
The permissibility conditions are $\e > -2$ and $\km \in \Rr_{+}$.

\subsection{Lommel Functions}
The \emph{Lommel functions} $S_{\mu,\nu}(z)$ are needed to solve the spectral integral~\eqref{eq:cov-BL-spectral}.
These functions are solutions of the inhomogeneous Bessel equation
\[
z^2 \, \frac{d^{2} w(z)}{dz^2} + z \frac{d w(z)}{dz} + (z^2 - \nu^2) \, w(z) = z^{\mu+1}.
\]
If either  $\mu + \nu$ or $\mu - \nu$ are an odd positive integer,  the respective Lommel functions are expressed as a terminating series, which
is given in descending order in the powers  of $z$ by the following equation~\cite[p.347]{Watson95}
\begin{align}
\label{eq:Smn-series}
S_{\mu,\nu}(z) &    = z^{\mu -1} \left[ \, 1 - \frac{(\mu-1)^2 - \nu^2}{z^2}  \right.
\nonumber \\
&   \left.  + \frac{[(\mu-1)^2 - \nu^2][(\mu-3)^2 - \nu^2]}{z^4} - \ldots \right]
\end{align}
If $\mu-\nu=2l+1$, where $l \in \mathbb{Z}_{+,0}$
the  series~\eqref{eq:Smn-series} terminates after $l+1$ terms.
In particular, the following general expression holds for such Lommel functions
\begin{align}
S_{\nu+2l+1,\nu}(z) &   =    z^{\nu+2l}   \,  \left( 1 + \mathbbm{1}_{l>0}(l) \sum_{k=1}^{l} (-1)^{k} \,
\right.
\nonumber \\
&   \cdot \left. \frac{\prod_{j=0}^{k-1} \left[(\nu +2(k-j))^2 - \nu^2 \right]}{z^{2k}}  \right).
\end{align}

\vspace{6pt}

\vspace{5pt}

\subsection{Bessel-Lommel Covariance in Position Space}{\hspace{16pt}}
\label{ssec:BL-cov}

\vspace{3pt}
\begin{proposition}
\label{prop:BL}
The Bessel-Lommel covariance
  $\inc^{BL}(z;\bmthe)$ as defined in~\eqref{eq:cov-BL-spectral} and~\eqref{eq:spd-BL}
is given by means of the following tripartite sum, where $z = \km \, \rr$,  $d\ge 2$, and $\nu=d/2-1$:

\begin{subequations}
\label{eq:cov-BL}
\begin{align}
\label{eq:Pz-all}
\inc^{BL}(z;\bmthe)  =  & \sum_{l=0,1,2} \frac{g_{l}(\bmthe)}{z^{2\nu+2l+1}} \,
\left[ (2\nu+2l) J_{\nu}(z)\, S_{\nu+2l,\nu-1}(z) \right.
\nonumber \\
&   \left.  \quad \quad  -  J_{\nu-1}(z)\, S_{\nu+2l+1,\nu}(z)  \right],
\end{align}
\begin{align}
\label{eq:g-l}
&   g_{0}(\bmthe) = \frac{ \km^{d}}{(2\pi)^{d/2} \eta _0\, \xi^{d} }, \;
g_{1}(\bmthe) = \e \, (\km \xi)^{2} \, g_{0}(\bmthe), \nonumber \\
&   g_{2}(\bmthe) = (\km \xi)^{4} \, g_{0}(\bmthe).
\end{align}
\end{subequations}
The  Lommel functions $S_{\mu,\nu}(z)$ used in~\eqref{eq:Pz-all} are shown in Table~\ref{tab:Lommel}.
For $ 0 < \km < \infty$ and  $\e >-2$, \eqref{eq:Pz-all} and~\eqref{eq:g-l}
define non-negative definite,
 infinitely differentiable radial functions.
\end{proposition}

\begin{IEEEproof}
The proof is given in the Appendix~\ref{App:proof-BL}.
\end{IEEEproof}
The Bessel-Lommel covariance functions can be viewed as a generalization of the cardinal sine covariance function,
$C_{\rm xx}(r) \propto \sin r/r$. The latter provides
a typical example of oscillatory decline of correlations due to finite spectral cutoff~\cite{Lantuejoul02,Yaglom87}.
\begin{table}[!t]
\renewcommand{\arraystretch}{1.3}
\caption{Lommel functions $S_{\nu+2l,\nu-1}(z)$ and $S_{\nu+2l+1,\nu}(z)$ for $l=0,1,2$ used in
$\inc^{BL}(z;\bmthe)$ given by~\eqref{eq:cov-BL}. The expressions are based on~\eqref{eq:Smn-series} where  $d\ge 2$ is the
space dimension and $\nu=d/2-1$.}
\label{tab:Lommel}
\centering
\begin{tabular}{|ll|ll}\hline
Notation & Lommel function  \\ \hline  \hline
 $S_{\nu,\nu-1}(z)$   &  $z^{\nu-1}$  \\
 $S_{\nu+2,\nu-1}(z)$   & $z^{\nu+1} \left( 1 - \frac{4\nu}{z^2} \right)$  \\
 $S_{\nu+4,\nu-1}(z)$   &  $z^{\nu+3} \left[ 1 - \frac{8(1+\nu)}{z^2} +  \frac{32\nu(1+\nu)}{z^4} \right]$ \\
 $S_{\nu+1,\nu}(z)$   &   $z^{\nu}$ \\
 $S_{\nu+3,\nu}(z)$   &  $z^{\nu+2} \left( 1 - \frac{4(1+\nu)}{z^2} \right)$    \\
 $S_{\nu+5,\nu}(z)$   &  $z^{\nu+4} \left[ 1 - \frac{8(\nu+2)}{z^2} +  \frac{32(\nu+1)(\nu+2)}{z^4} \right]$
 \\
 \hline
\end{tabular}
\end{table}

\subsection{Bessel-Lommel Variance and Correlation Function}
  Based on~\eqref{eq:invcovft} and~\eqref{eq:spd-BL}, the variance $\inc^{BL}(0;\bmthe)$
  is obtained by the following spectral integral
 \begin{align}
 \label{eq:BL-var}
 \inc^{BL}(0;\bmthe) &  =  \int_{\Rr^d} \frac{d\bfk}{(2\pi)^d} \,
 \frac{\mathbbm{1}_{\kk < \km }(\kk)}{\eta_{0} \xi^d}  \left( 1  +  \e\, \xi^2 \, \kk^2 + \xi^4 \, \kk^4\right)
 \nonumber \\
& =  \frac{2^{1-d}\, \km^d}{\pi^{d/2} \, \Gamma(d/2)\, \eta_{0} \, \xi^d} \left[ \frac{1}{d} +
 \frac{\e \, (\km \xi)^{2}}{d+2} +   \frac{(\km \xi)^{4}}{d+4} \right].
 \end{align}
Equation~\eqref{eq:BL-var} is obtained using the expression for the
 $d$-dimensional volume integral of radial functions $f(\kk)$, $\bfk \in \mathbb{R}^d$:
 \begin{equation*}
 \label{eq:iso-int}
 \int_{\Rr^d} d\bfk \, f(\kk) = \mathcal{S}_d \, \int_{0}^{\infty} dk \, k^{d-1} \, f(k),
 \end{equation*}
 where  $ \mathcal{S}_d = 2 \pi^{d/2}/\Gamma(d/2)$ is the surface of the unit sphere
 in $d$ dimensions~\cite[p. 39]{Schwartz08} and $\Gamma(\cdot)$ is the \emph{Gamma function}.

 Based on~\eqref{eq:BL-var},  $\inc^{BL}(0;\bmthe) \ll 1$ if $\eta_0\gg1$, or if
$\km /\sqrt{2\pi}\xi <1$ and  $d\gg1.$
The spatial dependence of the Bessel-Lommel covariance  is better understood
using the \emph{autocorrelation function}
$\rho_{\rm xx}^{BL}(z;\bmthe')$ ---where $\bmthe'=(\e, \xi, \km)^T $--- defined by
\begin{equation}
\label{eq:rho-BL}
\rho_{\rm x\rm x}^{BL}(z;\bmthe')= \frac{\inc^{BL}(z;\bmthe)}{\inc^{BL}(0;\bmthe)}.
\end{equation}

 \begin{figure}[!t]
  \centering
    \includegraphics[width=0.95\linewidth]{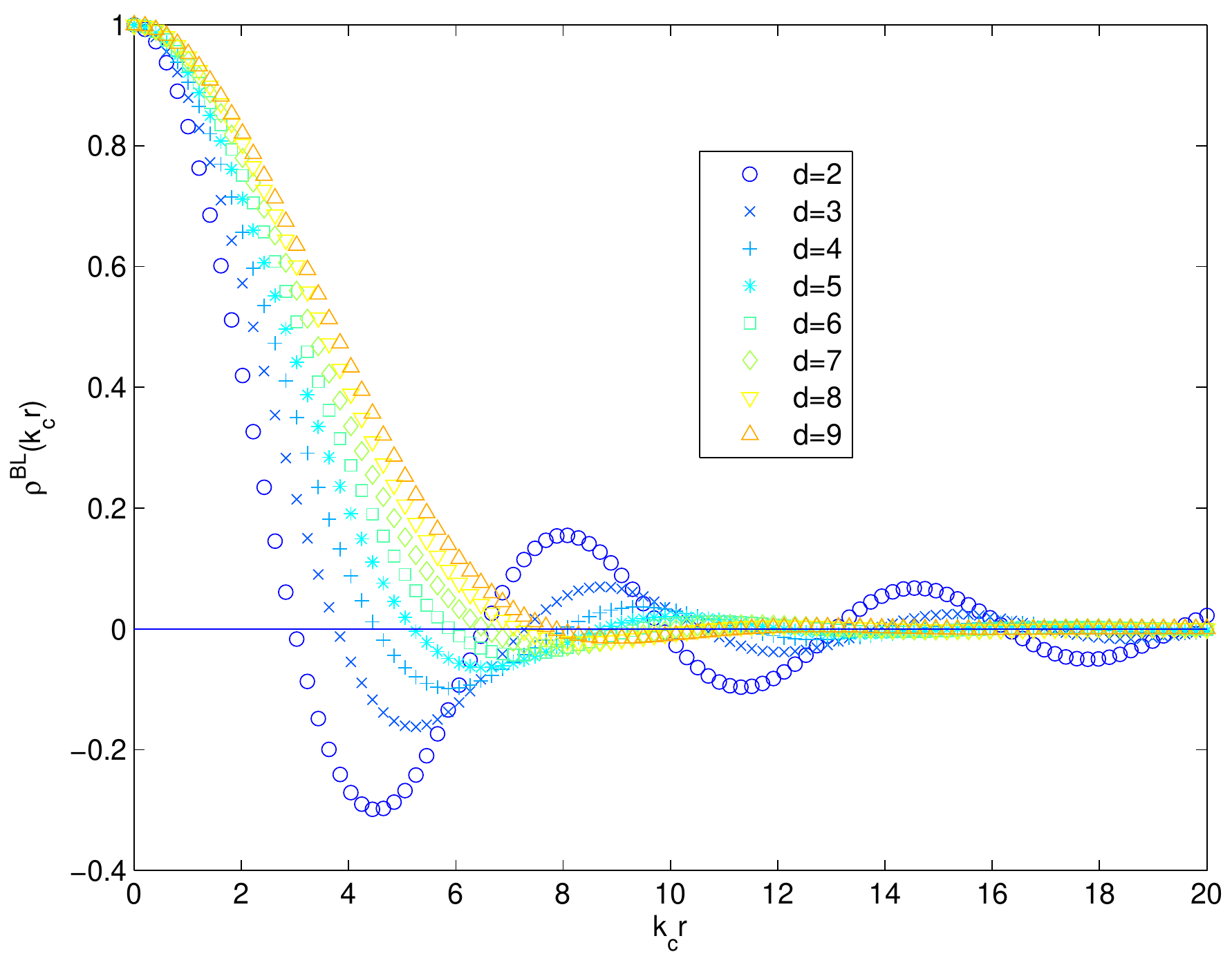}
  \caption{ Dependence of the Bessel-Lommel autocorrelation function~\eqref{eq:rho-BL} on the spatial lag and $d$.
  The values of the parameters are $\eta_0=1$, $\xi=1$, $\e=2$, and $\km=2$.}
  \label{fig:norm_invcov_ssrf_vs_d}
\end{figure}

The dependence  of $\rho_{\rm xx}^{BL}(z;\bmthe')$ on  $d$ is shown in
Fig.~\ref{fig:norm_invcov_ssrf_vs_d}. The amplitude of the negative peak of $\inc^{BL}(\rr;\bmthe)$
diminishes with increasing $d$, whereas the lowest positive root of $\inc^{BL}(\rr;\bmthe)$  moves to higher values.

\section{Length Scales of SSRF and Bessel-Lommel Random Fields}
\label{sec:cor-scales}
In two-parameter covariance models, e.g. the  Gaussian, $\Gxx(\rr)= \sigma^2 \, \exp(-\rr^2/\xi^2)$
and exponential,
 $\Gxx(\rr)= \sigma^2 \, \exp(-\rr/\xi)$ models, all the SRF length scales are determined by  $\xi$.
Models with more than two parameters, such as the SSRF, Mat\'{e}rn,
and  rational quadratic, exhibit more complex dependence since their length scales are not uniquely determined
by $\xi$~\cite{dth11}.

\begin{figure}[!t]
  \centering
    \includegraphics[width=0.85\linewidth]{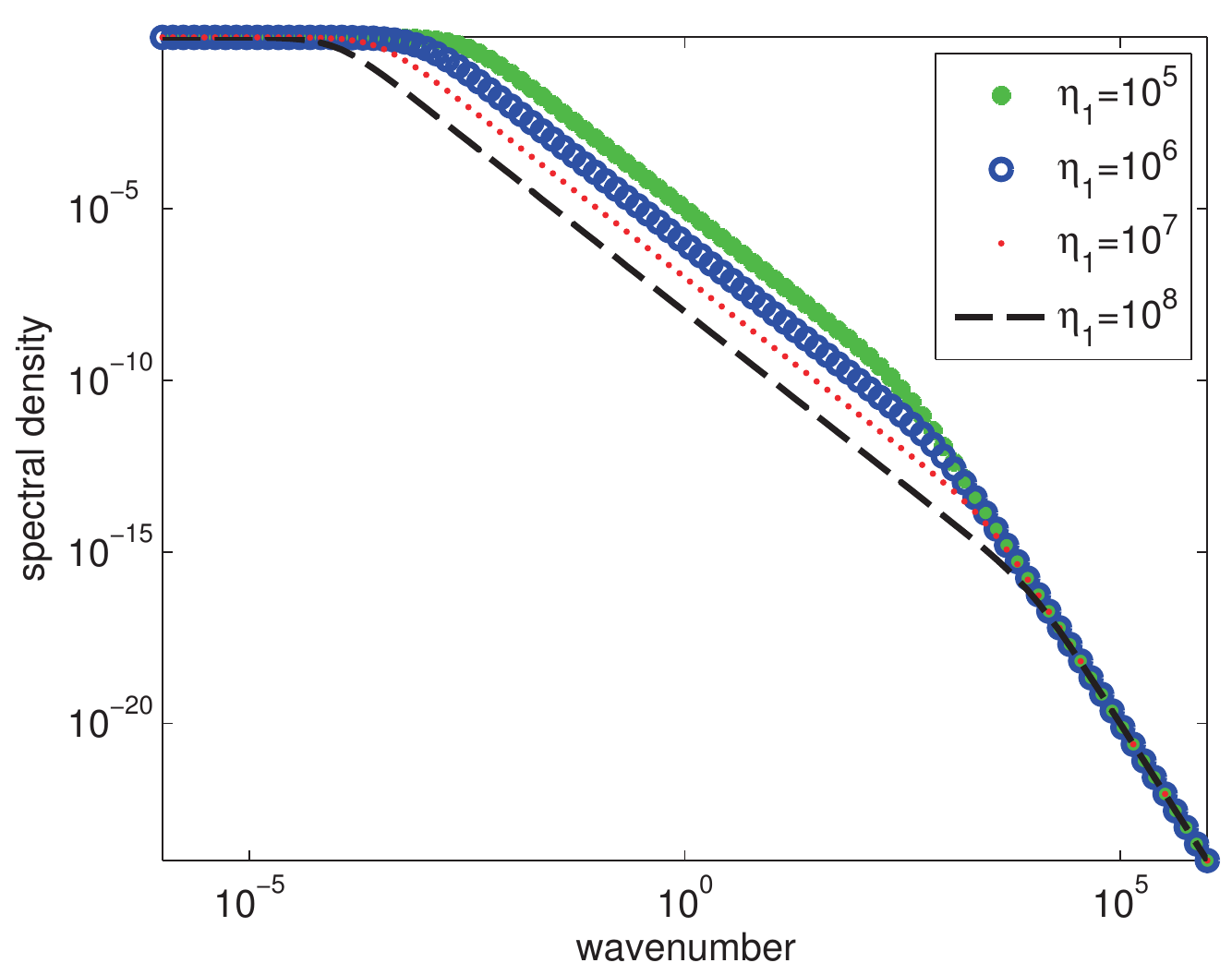}
  \caption{ Log-log plot of SSRF spectral density for $\eta_0=1$, $\xi=1$, and $\e=10^5, 10^6, 10^7, 10^8$.}
  \label{fig:ssrf-spd-mult}
\end{figure}

The SSRF spectral density exhibits power-law dependence over a wide range of
wavevectors if  $\e \gg 1$; in particular,
$\widetilde{\Gxx}(\bfk;\bmthe) \propto \kk^{-2}$  for $k_{\min} << \kk < k_{\max}$ where
$ k_{\min} = 1/(\sqrt{\e} \xi)$ and
$ k_{\max} = \sqrt{\e} /\xi$.  In Fig.~\ref{fig:ssrf-spd-mult} we show plots of the SSRF spectral
density for $\e \gg1$, which clearly exhibit the scaling (self-similar) range with an exponent equal to $-2$
over an extensive (but still finite) range.
Note that if the scaling dependence $\widetilde{\Gxx}(\bfk;\bmthe) \propto \kk^{-2}$
persisted for all $\kk \gg k_{\min}$ the
spectral density would not be summable.

In the Bessel-Lommel case, a necessary condition for self-similar scaling is $k_{\max} < \km.$
As we show below, however,  the
Bessel-Lommel integral range declines with $\km \uparrow$.  This implies that the wide spectral band
necessary to observe the scaling regime reduces spatial coherence.
Hence, the regime of self-similarity is not interesting in the Bessel-Lommel case.
\captionsetup[figure]{margin = 10pt}    

\subsection{Integral Range}
The \emph{integral range} represents a measure
of the distance over which two field values are correlated~\cite{Lantuejoul02,dth11}.
It is defined by the following volume integral
\begin{equation}
\label{eq:ellc}
 \ell_{c} = \left[  \frac{\int d{\bf r}\,
\Gxx(\rr ; \bmthe)}{\Gxx(0 ; \bmthe)} \right]^{1/d} =
\left( \frac{\widetilde{\Gxx}(\kk=0;\bmthe)}
{\int_{\Rr^d} d\bfk \,\widetilde{\Gxx}(\kk;\bmthe) }
\right)^{1/d}.
\end{equation}

\subsubsection{SSRF Covariance}
\label{ssec:ssrf-cov}
For SSRF covariance functions in $\Rr^2$ is $\ell_c = A_{2}(\e) \, \xi$ where $A_{2}(\e)$  is a monotonically increasing function of $\e$
 given by equations (20)-(22)  in~\cite{dth11}.

\subsubsection{Bessel-Lommel Covariance}

The integral range of the covariance function $\inc^{BL}(z;\bmthe)$ is given according to~\eqref{eq:ellc}
by the following expression
\begin{align}
\label{eq:BL-int-range}
 \ell_{c} = &   \left( \frac{\widetilde{\inc^{BL}}(\kk=0;\bmthe)}{\inc^{BL}(0;\bmthe)} \right)^{1/d}=
\frac{\pi^{1/2} \, 2^{1-1/d}}{ \km} \, [\Gamma(d/2)]^{1/d}
\nonumber \\
    &   \quad \cdot \left( \frac{1}{d} + \e \,
 \frac{(\km \xi)^{2}}{d+2} +   \frac{(\km \xi)^{4}}{d+4} \right)^{-1/d}.
 \end{align}
We derive~\eqref{eq:BL-int-range} using~\eqref{eq:spd-BL} for
$\widetilde{\inc^{BL}}(\kk=0;\bmthe)= (\eta_{0} \,\xi^d)^{-1}$ and~\eqref{eq:BL-var} for $\inc^{BL}(z=0;\bmthe).$

\begin{figure*}
\centering
\subfloat[Integral range vs. $d$]{\label{fig:intranged}
\includegraphics[angle=0, width=0.5\linewidth]{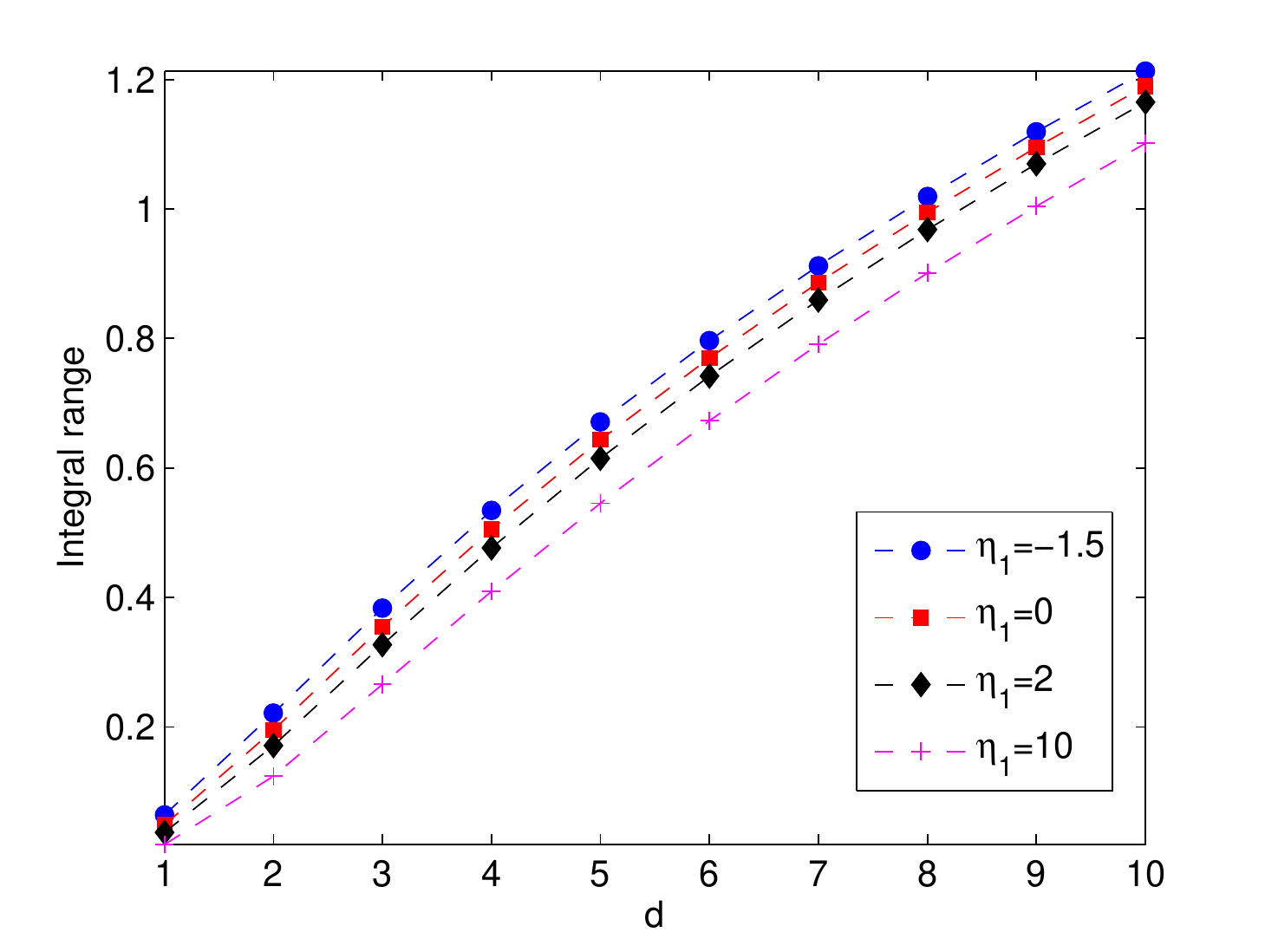}}
\subfloat[Integral range vs. $\km$]{\label{fig:intrangekc}
\includegraphics[angle=0, width=0.5\linewidth]{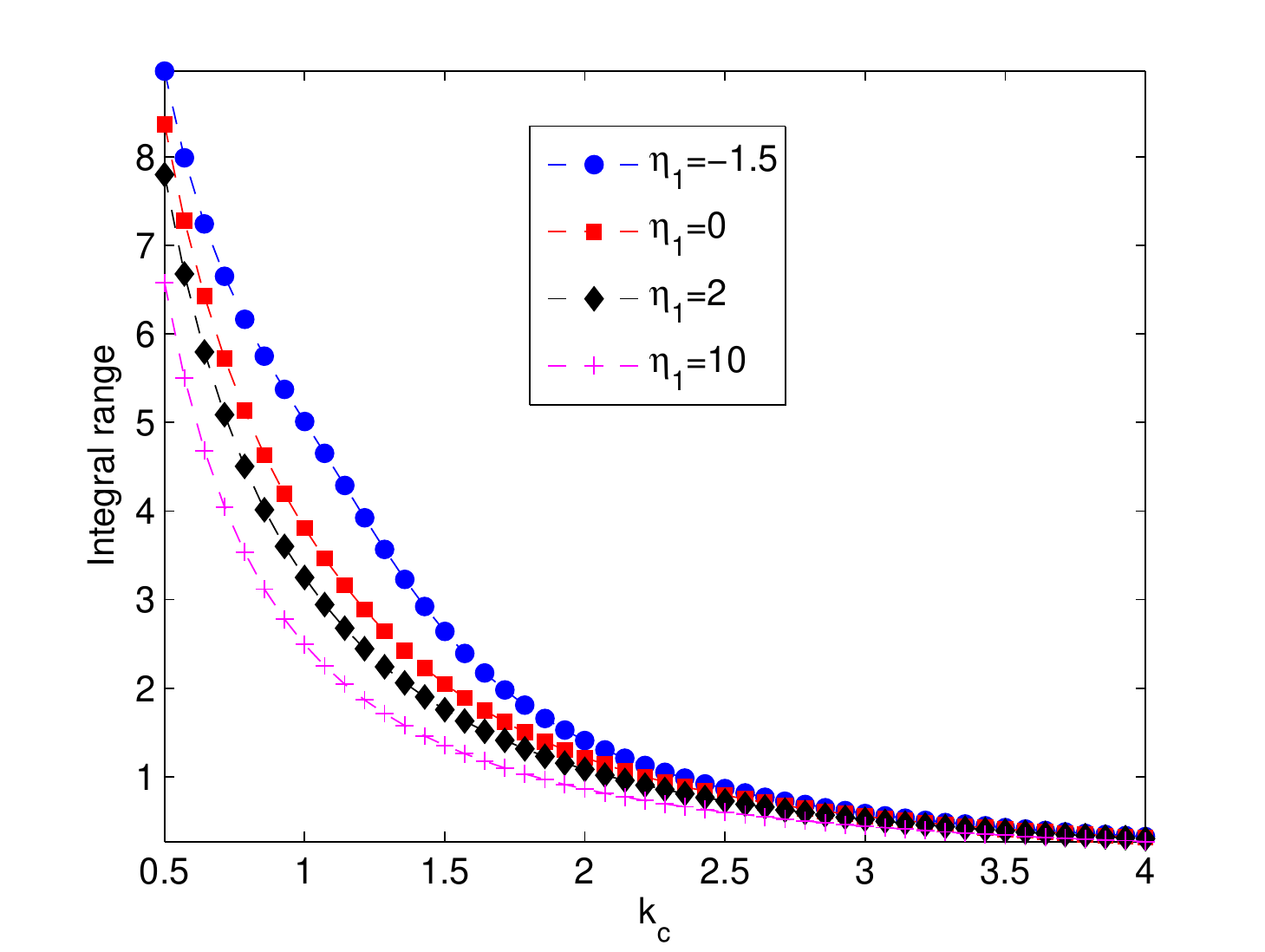}}
\caption{(a) Integral range $\ell_{c}$ versus  $d$
based on~\eqref{eq:BL-int-range}. (b) Integral range
 $\ell_{c}$ versus  $\km$  $\e$  based on~\eqref{eq:BL-int-range}; in both cases $\xi=2$.}
\end{figure*}

As evidenced in the parametric curves of Fig.~\ref{fig:intranged}, $\ell_c$ increases with $d$  reflecting the concomitant suppression of the oscillations
(see Fig.~\ref{fig:norm_invcov_ssrf_vs_d}).
In addition, $\ell_c$  decreases with increasing $\xi$ and $\e.$
The reason is the steeper increase of $\widetilde{\inc^{BL}}(\kk;\bmthe)$  with $\kk$
as $\xi \uparrow$ or $\e \uparrow$,
which enhances spatial variability and reduces the SRF coherence.
The above trend  is  opposite to that of the SSRF integral range.
The dependence of $\ell_c$ on $\xi$ enters~\eqref{eq:BL-int-range} only through the product $\km \xi= u_{c}$.
Hence, for fixed $u_c$, $d$ and $\e$, it follows from~\eqref{eq:BL-int-range} that  $ \ell_{c} \propto 1/\km $
as shown in~ Fig.~\ref{fig:intrangekc}.
The observed decline occurs  because  $\km \uparrow$ signifies a wider spectral band with increasing weight at its tail
which reduces spatial coherence.
The relation $ \ell_{c} \propto 1/\km $ is justified by the fact that
 $\inc^{BL}(z;\bmthe)$ is a function of $z= u_c \, h =  \km \rr$, which implies
that the characteristic distance scale is $\propto 1/\km$ instead of $\xi.$
\captionsetup[subfigure]{margin = 10pt}    

\subsection{Correlation Spectrum}

The reference~\cite{Kramer07}
introduces the concept of the
\emph{smoothness microscale}. The microscale denotes a length
such that the SRF appears smooth and  can be linearly interpolated at
smaller length scales. The microscale is equivalent to the integral range of
the gradient of $X(\bfs,\om)$.
We extend the definition in~\cite{Kramer07} to
a correlation spectrum  which applies to  stationary,  but not necessarily mean-square
differentiable SRFs with unimodal, radially symmetric  spectral density.

\vspace{6pt}
\begin{definition}
\label{defi:lac}
Let the radial function $\Rr^d \supset \bfk \mapsto \widetilde{\Gxx}(\kk;\bmthe) $
be a permissible spectral density for  a statistically isotropic SRF $X(\bfs,\om)$
that satisfies Bochner's theorem~\ref{theor:Bochner}.
 In addition, let
  $\widetilde{\Gxx}(\kk;\bmthe)$ be a unimodal function of $\kk$.
  Then, we define the following correlation spectrum $\lambda_{c}^{(\alpha)}$
  indexed by $ 0 \le \alpha \le 1$
\begin{align}
\label{eq:micro} \lambda_{c}^{(\alpha)} = &    \left(   \frac{\sup_{\bfk \in \Rr^d} \, \kk^{2\alpha} \, \widetilde{\Gxx}(\kk;\bmthe)}
{\int_{\Rr^d} d\bfk \,\kk^{2\alpha} \,\widetilde{\Gxx}(\kk;\bmthe) }  \;
\right)^{1/d}.
\end{align}
\vspace{6pt}

\end{definition}

For $\alpha = 0$ the \emph{integral range}~\eqref{eq:ellc} is recovered if the peak
of the spectral density occurs
at $\kk=0.$  If the spectral density reaches a maximum at wavenumber $k_{\max}> 0$, then~\eqref{eq:micro}
 represents approximately the width of the maximum of the spectral density.

For $\alpha =1 $ the \emph{smoothness microscale} is recovered.
The definition of the microscale involves in the denominator of~\eqref{eq:micro} the Laplacian of the
covariance function evaluated at zero lag, $-[\nabla^{2} \Gxx(\rr)]_{\rr=0}$.
If the SRF is differentiable in the mean square sense
$-[\nabla^{2} \Gxx(\rr)]_{\rr=0} \in \Rr_{+}$,  and thus
$\lambda_{c}^{(\alpha=1)} \in \Rr_{+}$;
in contrast, $[\nabla^{2} \Gxx(\rr)]_{\rr=0}$  diverges for non-differentiable SRFs leading to $\lambda_{c}^{(\alpha=1)}=0$. The zero  microscale denotes that
the SRF appears rough at all length scales.

Exponents $0 < \alpha <1$ generate a spectrum of
scales which emphasizes different parts of the  spectral density.  Correlation  lengths obtained
 for $0< \alpha <1$ correspond to the integral range of the SRF's fractional derivative of order $\alpha$.
Based on~\eqref{eq:micro} these length scales involve the fractional Laplacian of the covariance function~\cite{Samko93} at zero lag,
and they take finite, non-zero values
if $\kk^{2\alpha} \,\widetilde{\Gxx}(\kk;\bmthe) $ is integrable. A necessary condition for integrability
is that $\widetilde{\Gxx}(\kk;\bmthe) \underset{\kk \ra \infty}{\sim} \kk^{-q}$  where $ q > 2\alpha + d $.
Hence, for fixed $q$ a non-vanishing length scale is obtained for $\alpha < \alpha_{\max} = (q - d)/2 $;
thus, the respective length scales   $\lambda_{c}^{(\alpha)}$ can be used to quantify
the ``fractional smoothness'' of continuous but non-differentiable random fields.

\subsubsection{SSRF Covariance}
For $\alpha=1$, it follows from~\eqref{eq:micro} that $\lambda_{c}^{(\alpha)}=0$, because the integral
in the denominator develops a logarithmic divergence
marking the non-existence of the  covariance Laplacian at zero.
The loss of mean square differentiability implies that
the smoothness microscale is zero
regardless of $\e$ and $\xi$.  The failure of the microscale to discriminate  between different
SSRF parameters is unsatisfactory, since the latter
affect the small-scale structure of field realizations as evidenced in the
plots shown in Fig.~\ref{fig:ssrf-simul} below.

For $0 \le \alpha <1$, however, $\lambda_{c}^{(\alpha)}> 0$; in addition,
$\lambda_{c}^{(\alpha)}$ depends on the SSRF parameters.
 Hence, $\lambda_{c}^{(1 -\epsilon)}$ where $ 0 < \epsilon \ll 1$  distinguishes between different SSRF covariance models.

 \begin{proposition}
 \label{prop:ssrf-scales}
 The  correlation spectrum for the SSRF covariance functions with spectral density~\eqref{eq:ssrf-spd}
 are given by the following equation in $3 \ge d\ge 1$
\begin{align}
\label{eq:micro-ssrf}
\lambda_{c}^{(\alpha)} = &  \xi \, \left(  \frac{(\tilde{\kappa}_{1}\, \xi)^{2\alpha} \Big/
b(\alpha,\e)}{1 + \e (\tilde{\kappa}_{1}\xi)^2
 + (\tilde{\kappa}_{1}\xi)^4 }\right)^{1/d},
 \end{align}
 where
\begin{align}
\label{eq:micro-ssrf-2}
b(\alpha,\e) &  =
 \left\{ \begin{array}{cc}
           \frac{\pi \,\Gamma(1-\alpha)\Gamma(1+\alpha) \, \left[\left( \e + \Delta \right)^{\alpha}  -
\left( \e - \Delta\right)^{\alpha}\right]}{2^{\alpha}\alpha \Delta}    &   \e \neq 2  \\
           \pi \, \Gamma(1+\alpha) \, \Gamma(1-\alpha)   &   \e=2
        \end{array}  \right.
\end{align}
$\Delta = \sqrt{\e^2 -4}$, and $\tilde{\kappa}_{1} = \underset{{\kk}}{\arg \max} \left( \kk^{2\alpha} \, \widetilde{\Gxx}(\kk;\bmthe) \right)$, i.e.,
\begin{equation}
\label{eq:root-fk}
\tilde{\kappa}_{1} = \sqrt{\frac{ \sqrt{\e^2 \, (1 - \alpha)^2 - 4 \alpha (\alpha-2)} -
\e \, (1 - \alpha)}{2(2 - \alpha) \xi^2}}.
\end{equation}
\end{proposition}

\begin{IEEEproof}
The proof is given in the Appendix~\ref{App:Proof-C}.
\end{IEEEproof}

\begin{figure}[!t]
\centering
\subfloat[SSRF scales - fixed $\e$]{\label{fig:ssrf-cora-alpha}
\includegraphics[angle=0, width=0.5\linewidth]{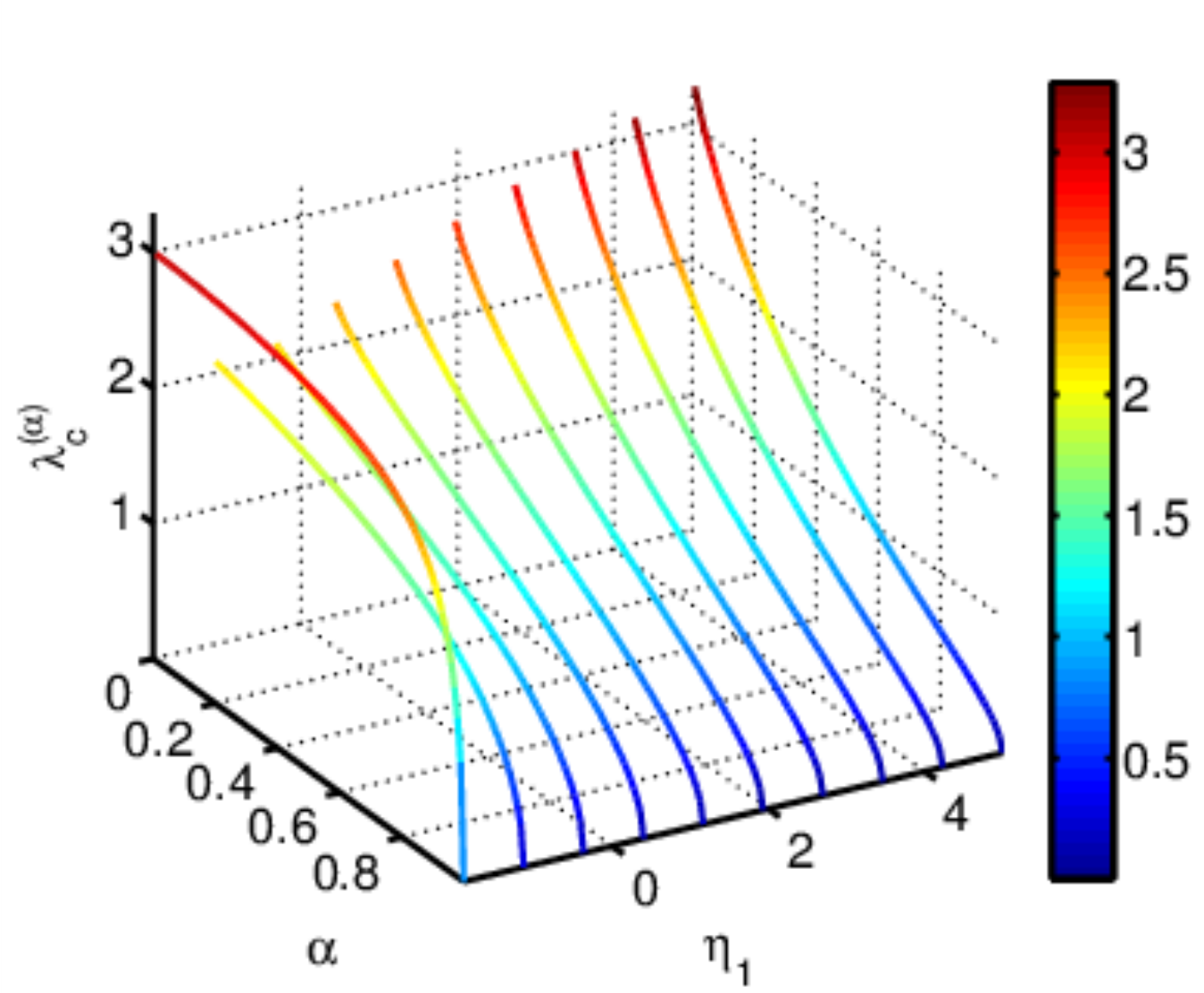}}
\subfloat[SSRF scales - fixed $\alpha$]{\label{fig:ssrf-cora-eta1}
\includegraphics[angle=0, width=0.5\linewidth]{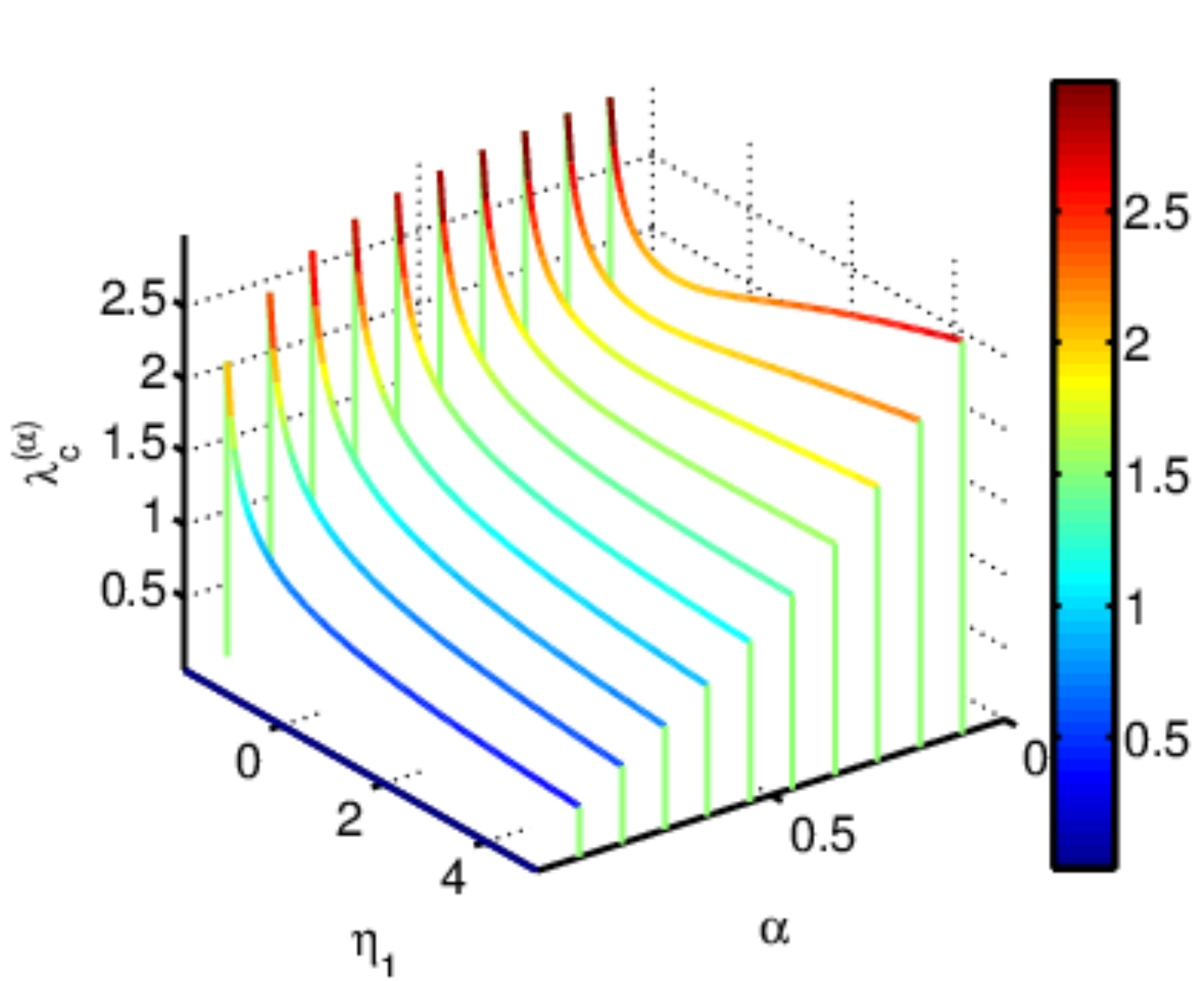}}
\caption{ SSRF correlation spectrum~\eqref{eq:micro-ssrf} versus
$\e \in [-1.9, 5]$ and $\alpha \in [0, 1]$
for  $\xi=5$ and $d=2$. (a) Continuous lines correspond to curves of fixed $\e$.
The curves displayed are for ten evenly spaced $\e$ values between $-1.9$ to $5$ (from left to right in the plot).
(b) Continuous lines correspond to curves of fixed $\alpha$. The curves displayed are for ten evenly spaced $\alpha$ values between
$0.05$ to $0.95$ (from right to left in the plot).
 }
\label{fig:SSRF-spectrum-scales}
\end{figure}

The dependence of $\lambda_{c}^{(\alpha)}$ on $\e$  and $\alpha$ is shown in Fig.~\ref{fig:SSRF-spectrum-scales}.
The general trend, as shown in Fig.~\ref{fig:ssrf-cora-alpha}, is that $\lambda_{c}^{(\alpha)} \downarrow$  as $\alpha \uparrow$, leading to $\lambda_{c}^{(1)}=0$ which marks the
SSRF non-differentiability. The decline of $\lambda_{c}^{(\alpha)}$ with $\alpha \uparrow$ reflects the fact that
higher $\alpha $ emphasize shorter length scales (wavelenghts) as marked by the increase of $\tilde{\kappa}_{1}$ with $\alpha$ at fixed $\e$.
The curves in Fig.~\ref{fig:ssrf-cora-eta1} show that for  $\alpha > \alpha_{\min} \approx 2$,
$\lambda_{c}^{(\alpha)} \downarrow$  as $\e \uparrow$,
indicating that higher rigidity reduces the fine-scale regularity of the field.
For $\alpha \le \alpha_{\min} $, i.e., for coarser scales,
first $\lambda_{c}^{(\alpha)} \downarrow$  for $\e  \lessapprox -1$, and then
 $\lambda_{c}^{(\alpha)} \uparrow$ as $\e \uparrow$  beyond $\e \gtrapprox -1$.

\subsubsection{Bessel-Lommel Covariance}
\label{sssec:cor-scales-BL}

\begin{proposition}
 \label{prop:BL-scales}
 The correlation spectrum for the Bessel-Lommel covariance functions with  spectral density~\eqref{eq:spd-BL} are
 given by the following equation in $d\ge 2$
\begin{subequations}
\label{eq:BL-micro}
\renewcommand\minalignsep{3pt}
\begin{flalign}
\lambda_{c}^{(\alpha)}  & =
\frac{g_{d}}{\km } \left(   1 + \e \tilde{\km} ^2 +
\tilde{\km}^4 \right)^{1/d}, \;   \e> \ec  \lor  \km < \kappa_{-}
\label{eq:BL-micro-pos-eta1}
\\
\lambda_{c}^{(\alpha)}  & = \frac{g_{d}}{\km } \left(   1 + \e \tilde{\kappa}_{-}^2 +
\tilde{\kappa}_{-}^4 \right)^{1/d}, \;  \ec \ge \e> - 2
\nonumber \\
&   \hspace{1.5in}     \wedge  \kappa_{+} \ge \km > \kappa_{-},
\label{eq:BL-micro-neg-eta1-1}
\\
\lambda_{c}^{(\alpha)}  &  = \frac{g_{d}}{\km } \left(   1 + \e (\tilde{\kappa}^{\ast} )^2 +
(\tilde{\kappa}^{\ast} )^4 \right)^{1/d}, \;  \ec \ge \e> - 2  \nonumber \\
&   \hspace{1.6in}  \wedge  \km > \kappa_{+},
\label{eq:BL-micro-neg-eta1-2}
\end{flalign}
\end{subequations}
where $\tilde{\km} = \km \xi$, $\tilde{\kappa}_{-}= {\kappa}_{-} \xi$, $\tilde{\kappa}^{\ast}={\kappa}^{\ast}\xi$,
\begin{subequations}
\begin{equation}
\ec = - \frac{2\,\sqrt{\alpha(\alpha+2)}}{(\alpha+1)},
\end{equation}
\begin{equation}
g_{d} = \left[ \frac{\Gamma(d/2)}{2\, \pi^{d/2} \,
\left( \frac{1}{d+2\alpha} + \frac{\e \,\km^2 \,\xi^2}{d+2\alpha+2} + \frac{\km^4 \,\xi^4}{d+2\alpha+4} \right)}
\right]^{1/d},
\end{equation}
\begin{equation}
\label{eq:BL-kappa}
\kappa_{\pm} = \frac{1}{\xi}\sqrt{\frac{-(\alpha+1)\e \pm \sqrt{(\alpha+1)^2 \e^2 - 4\alpha(\alpha+2)}}{2 (\alpha+2)}},
\end{equation}
\end{subequations}
and $\kappa^{\ast}= \arg \max \left( \kappa_{-}^{2\alpha}\, \Pi(\kappa_{-} \xi),  \km^{2\alpha}\, \Pi(\km \xi) \right)$, $\Pi(x)$ being the polynomial defined in~\eqref{eq:char-pol}.
\end{proposition}

\begin{IEEEproof}
 The proof is given in the Appendix~\ref{App:Proof-D}.
\end{IEEEproof}

\begin{figure}[!t]
\centering
\subfloat[BL scales $\xi=5, \km=\pi/2$]{\label{fig:bl-cora-e1}
\includegraphics[angle=0, width=0.5\linewidth]{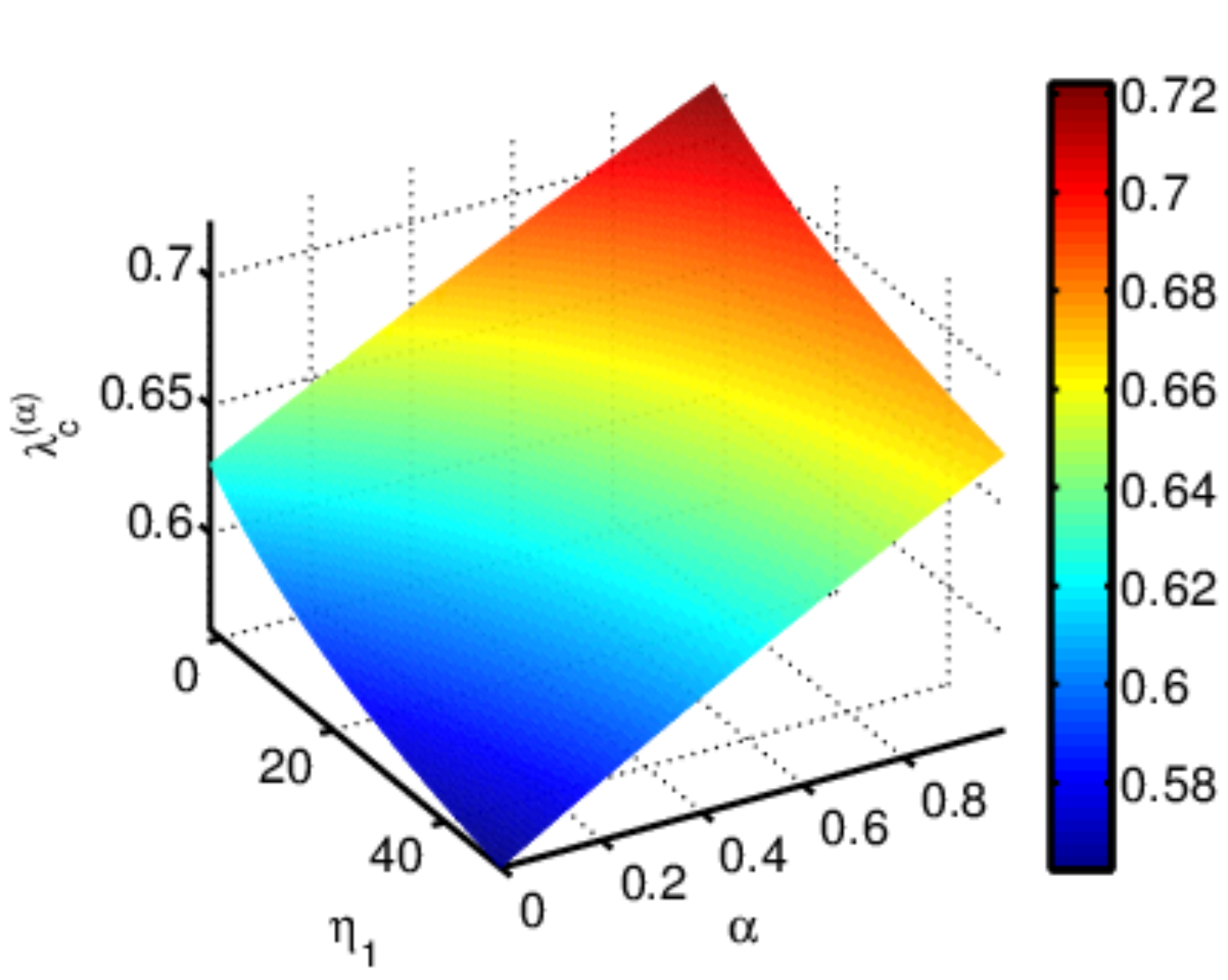}}
\subfloat[BL scales $\e=3, \xi=5$]{\label{fig:bl-cora-kc}
\includegraphics[angle=0, width=0.5\linewidth]{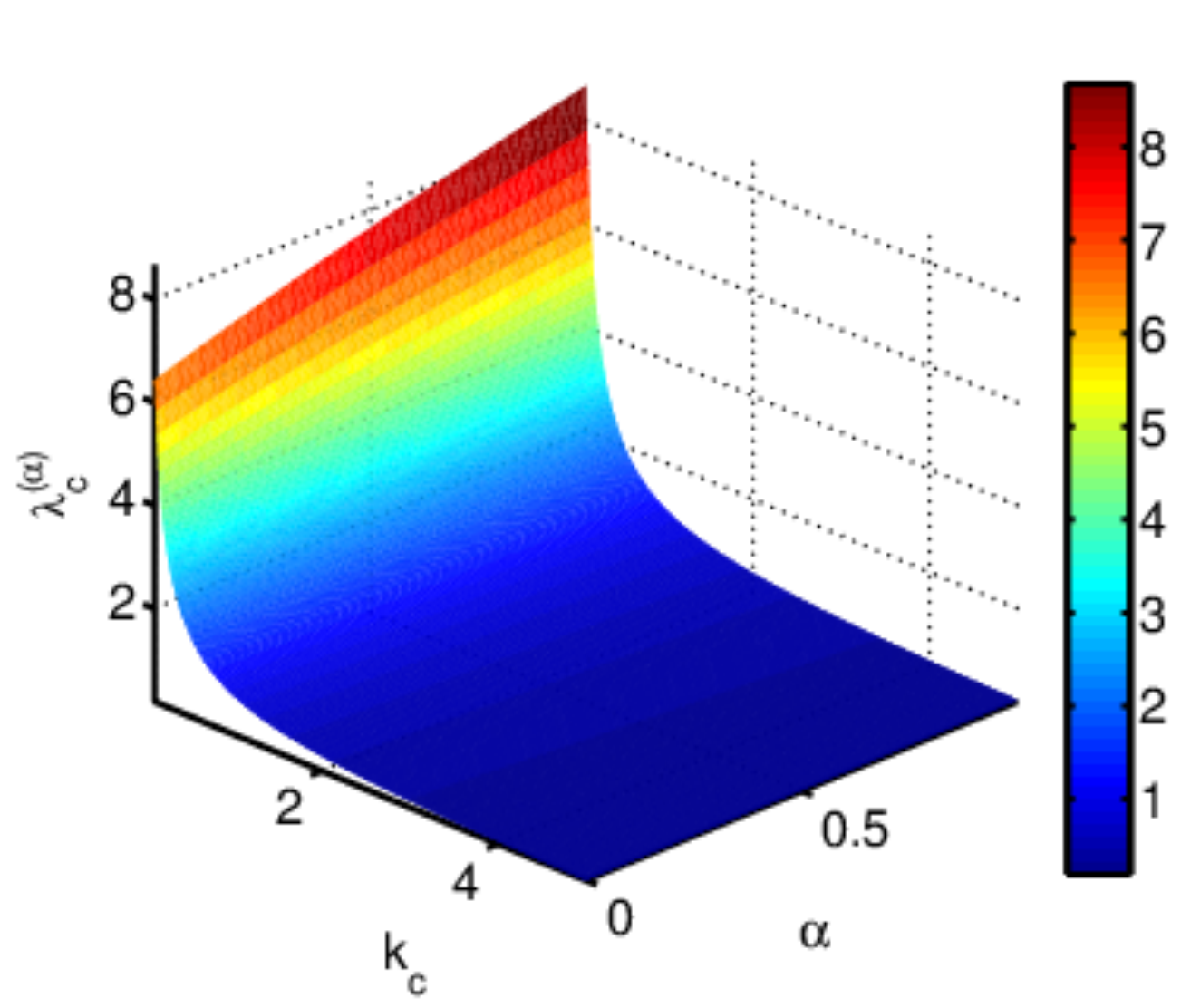}}
\caption{Spectrum of Bessel-Lommel correlation scales in $d=2$. (a) $\lambda_{c}^{(\alpha)}$ for
$\e \in [0, 50]$ and $\alpha \in [0, 1]$
with  $\xi=5, \km=\pi/2$. (b) $\lambda_{c}^{(\alpha)}$ for
$\km \in [0.1, 5]$ and $\alpha \in [0, 1]$
with  $\e=3, \xi=5$.}
\label{fig:BL-spectrum-scales}
\end{figure}

The dependence of $\lambda_{c}^{(\alpha)}$ on $\e, \alpha, \km$ is shown in Fig.~\ref{fig:BL-spectrum-scales}.
Overall, higher $\e$ and $\km$ tend to reduce $\lambda_{c}^{(\alpha)}$. The $\km$ dependence is more pronounced and follows asymptotically
$\lambda_{c}^{(\alpha)} \sim 1/\km^{3}$. The increase of $\lambda_{c}^{(\alpha)}$
with $\alpha$ signifies that
the smoothness microscale exceeds the integral range. This effect is due
on one hand to the high smoothness of the B-L random field that increases the microscale and
on the other to the oscillations of the B-L covariance function which reduce the integral range.

\section{Simulations}
\label{sec:simul}
We present realizations of SRFs with SSRF and Bessel-Lommel correlation structures.
The spectral simulation method based on the Fast Fourier Transform (FFT)~\cite{Gutjahr93,Kramer07,dth05} is used to
generate the realizations on square grids of length $L=512$.

\subsection{Random Fields with SSRF Covariance}

\begin{figure*}
\centering
\subfloat[$\eta_{1}=-1.5$]{\label{fig:ssrf_1}
\includegraphics[angle=0, width=0.25\linewidth]{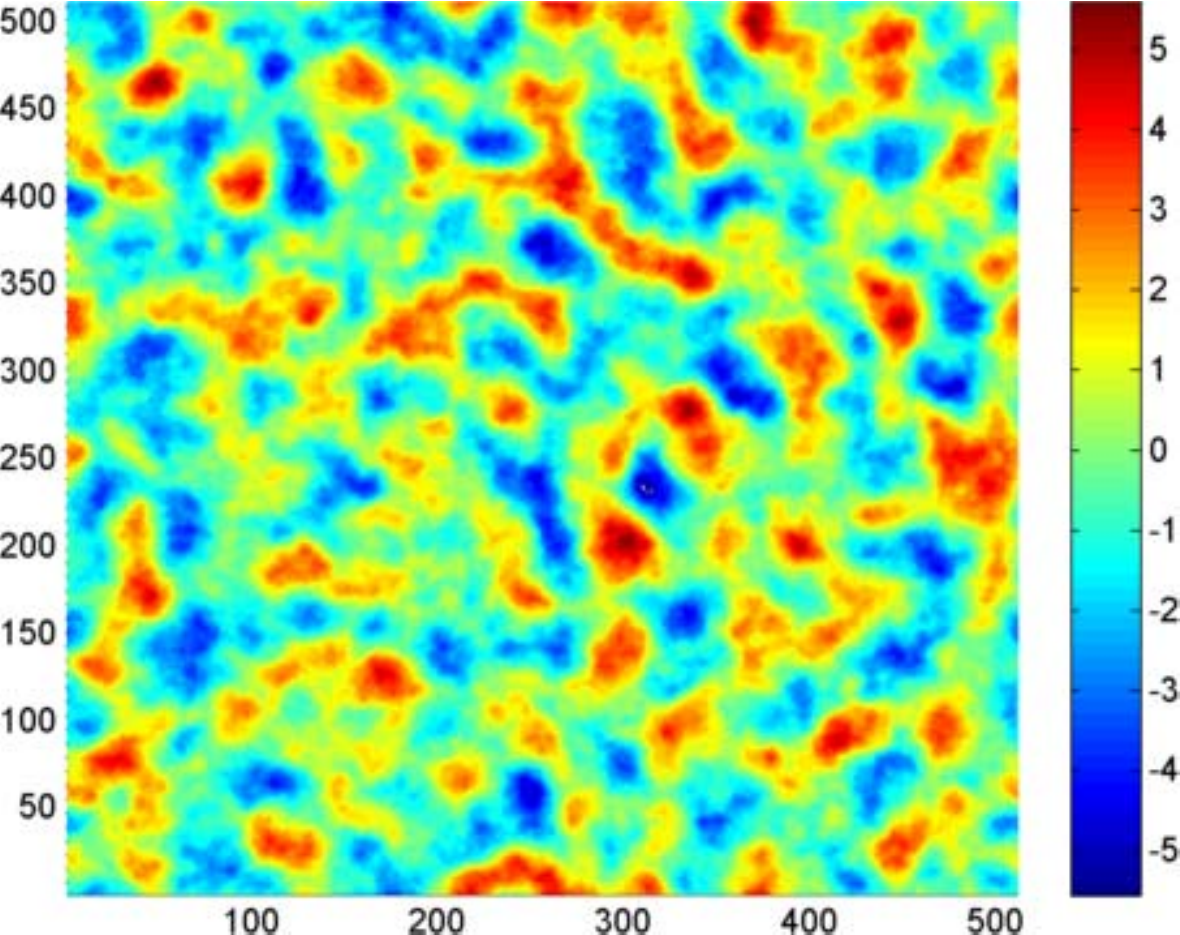}}
\subfloat[$\eta_{1}=0$]{\label{fig:ssrf_2}
\includegraphics[angle=0, width=0.25\linewidth]{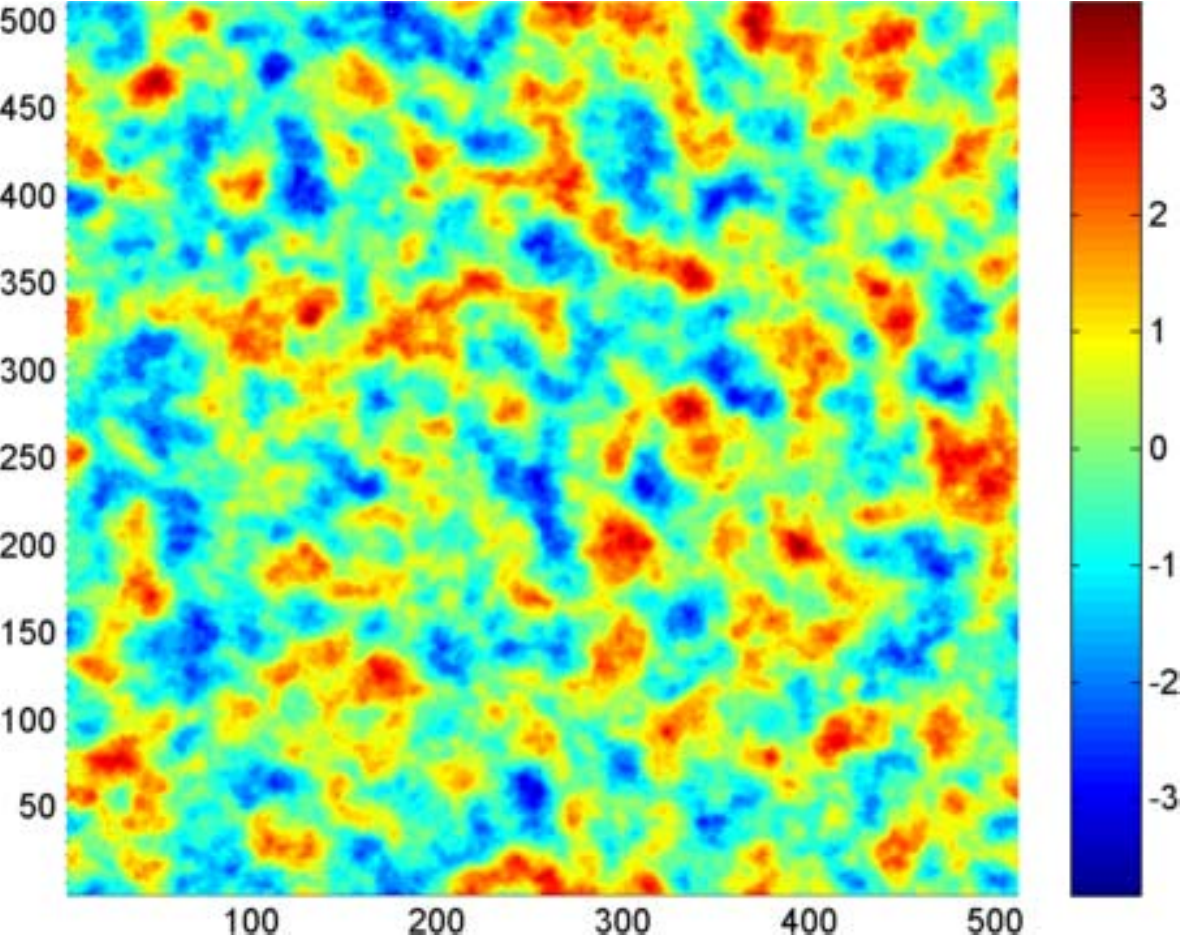}}
\subfloat[$\eta_{1}=1.5$]{\label{fig:ssrf_3}
\includegraphics[angle=0, width=0.25\linewidth]{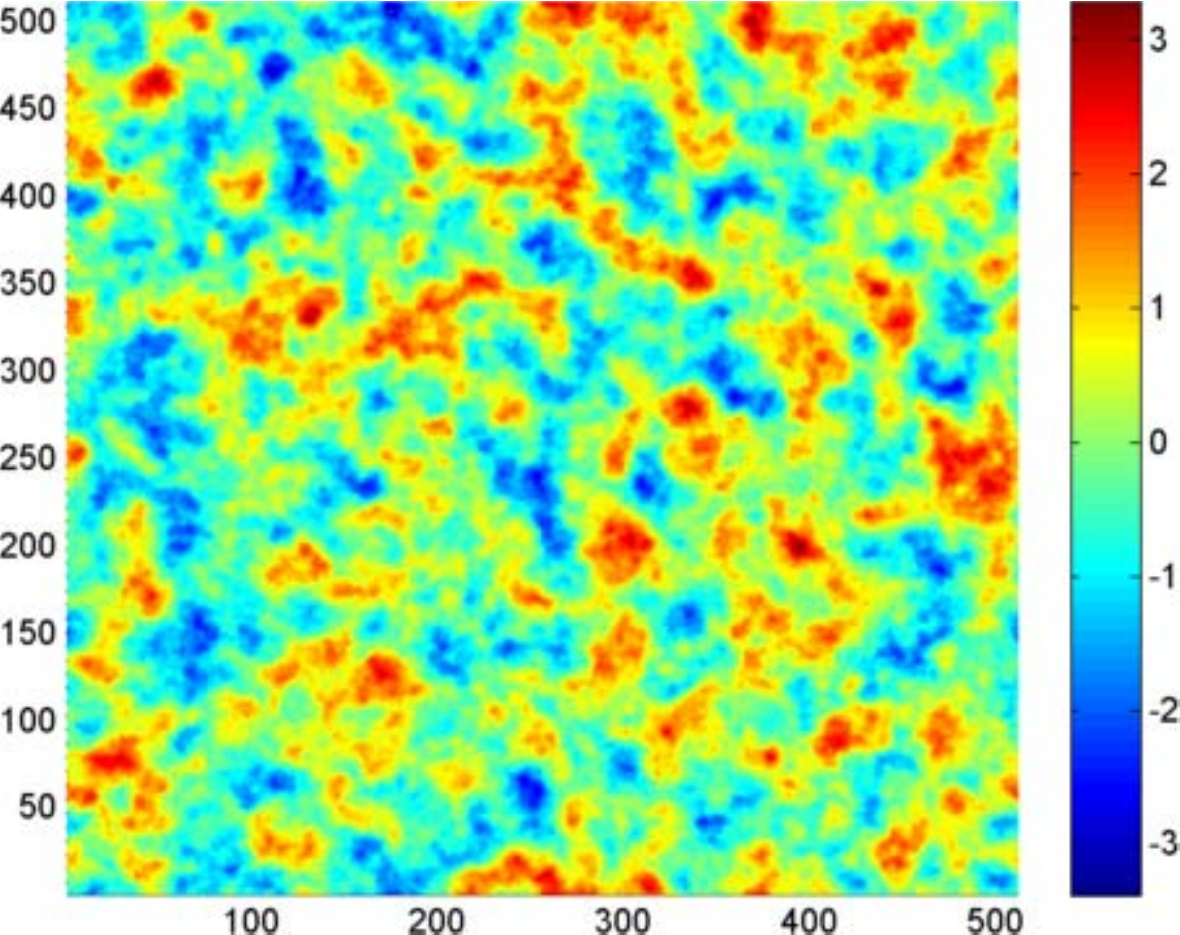}}
\subfloat[$\eta_{1}=15$]{\label{fig:ssrf_4}
\includegraphics[angle=0, width=0.25\linewidth]{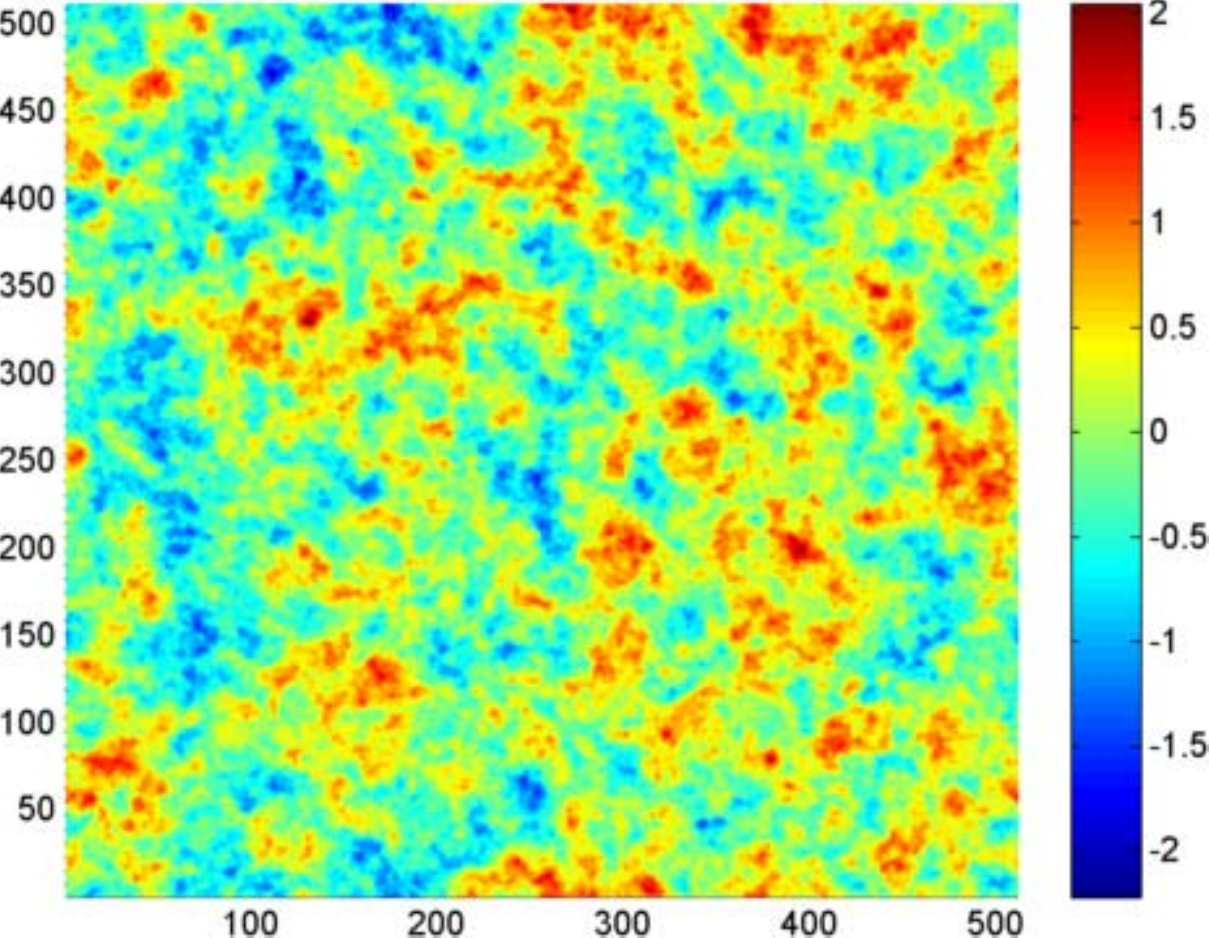}}
\caption{SSRF realizations with covariance function~\eqref{eq:cov2d-ssrf} on a square $512 \times 512$ grid.  For all  realizations  $\eta_{0}=10, \xi=10$.
 FFT spectral simulation  with the same random generator seed is used in all the simulations. }
 \label{fig:ssrf-simul}
\end{figure*}

Four different SSRF realizations
are shown in Fig.~\ref{fig:ssrf-simul}. The realizations correspond to
fields with common $\eta_0$ and $\xi$ but different $\eta_1$.
The patterns become ``grainier'' as $\eta_{1}$ increases. Nevertheless, the integral range $\ell_c$
increases with $\eta_1$,  as discussed in section~\ref{ssec:ssrf-cov}.
\captionsetup[subfigure]{margin = 10pt}    

\begin{figure}[!t]
\centering
\subfloat[$\eta_{0}=10^4, \e=10^5, \xi=10$]{\label{fig:ssrf_n_1}
\includegraphics[angle=0, width=0.5\linewidth]{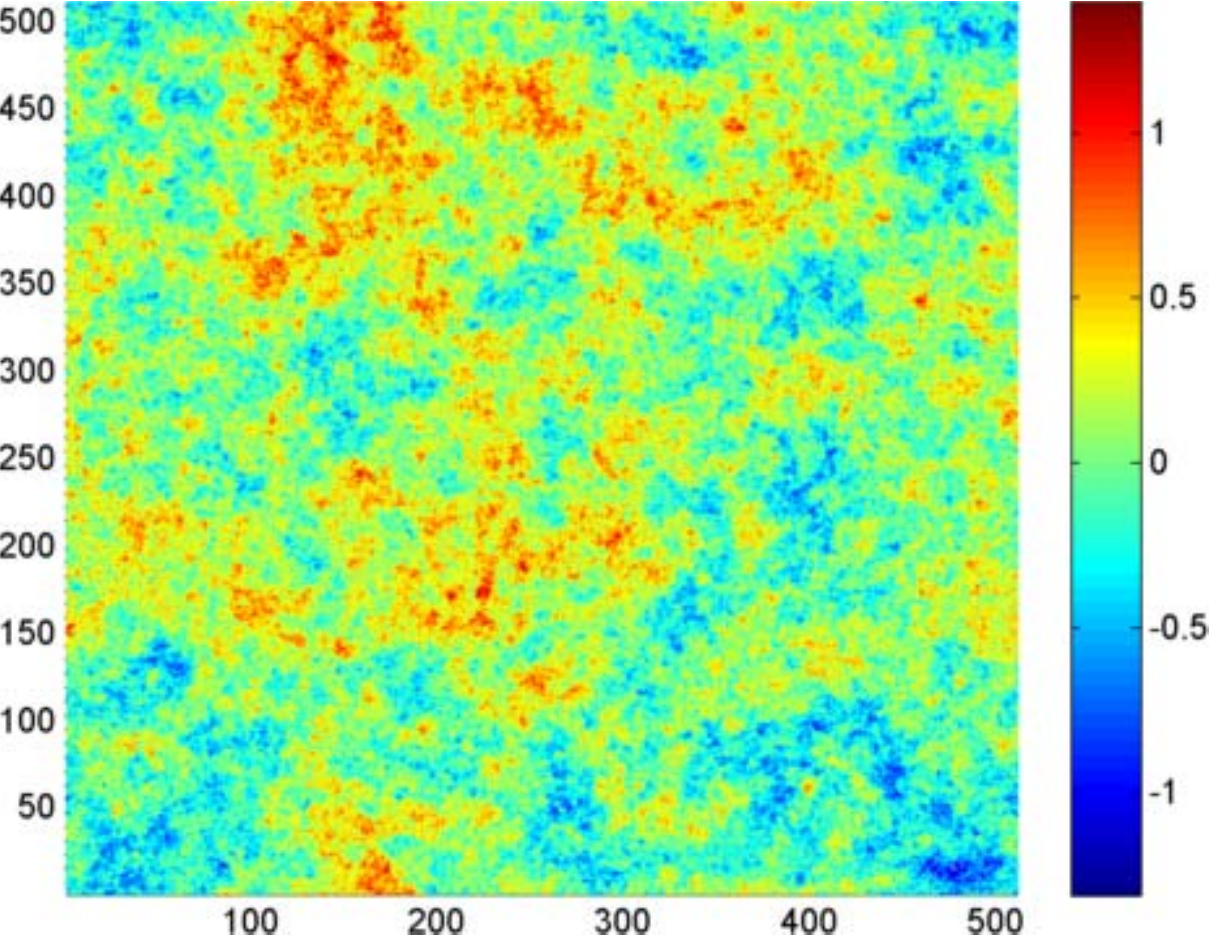}}
\subfloat[$\eta_{0}=10^6, \e=10^6, \xi=10$]{\label{fig:ssrf_n_2}
\includegraphics[angle=0, width=0.5\linewidth]{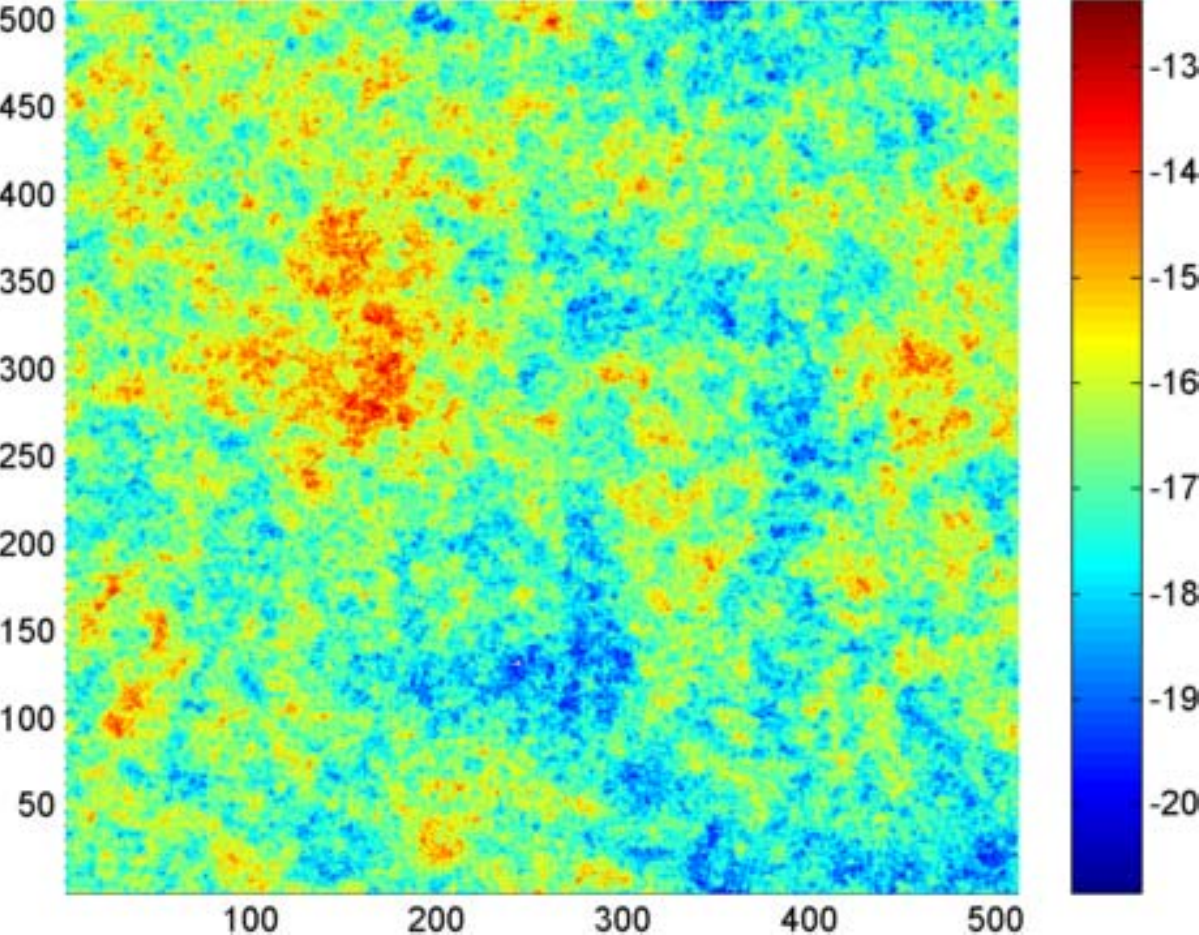}}
\caption{Multiscale SSRF realizations with covariance function~\eqref{eq:cov2d-ssrf} on square $512 \times 512$ grid.  (a):
$\eta_{0}=10^4, \e=10^5, \xi=10$; (b): $\eta_{0}=10^6, \e=10^6, \xi=10$.
FFT spectral simulation with clock dependent random generator seed is used.  }
 \label{fig:ssrf-simul-multi}
\end{figure}

Two realizations with $\e \gg 1$ that exhibit multiscaling  are shown in Fig.~\ref{fig:ssrf-simul-multi}.
For $\e \gg 1$ it follows from~\cite[eq. (22)]{dth11} that $r_c \approx  \xi\sqrt{2\pi \, \e/ \log \e}$.
Hence, the realization shown in Fig.~\ref{fig:ssrf_n_1}
is drawn from a field with $r_c \approx 2336$, whereas the one shown in Fig.~\ref{fig:ssrf_n_2} corresponds to $r_c \approx 6743$. In both cases the domain size does not allow adequate
sampling of the probability distribution. This non-ergodic effect is
 responsible for the difference
in the range of field values between the two plots.
The lack of ergodicity is more pronounced in Fig.~\ref{fig:ssrf_n_2},
where all the field values  are  negative (see vertical scale bar).
\captionsetup[subfigure]{margin = 10pt}    

\subsection{Random Fields with Bessel-Lommel Covariance}
\begin{figure*}
\centering
\subfloat[$\km=0.05$]{\label{fig:BesLom_1}
\includegraphics[angle=0, width=0.25\linewidth]{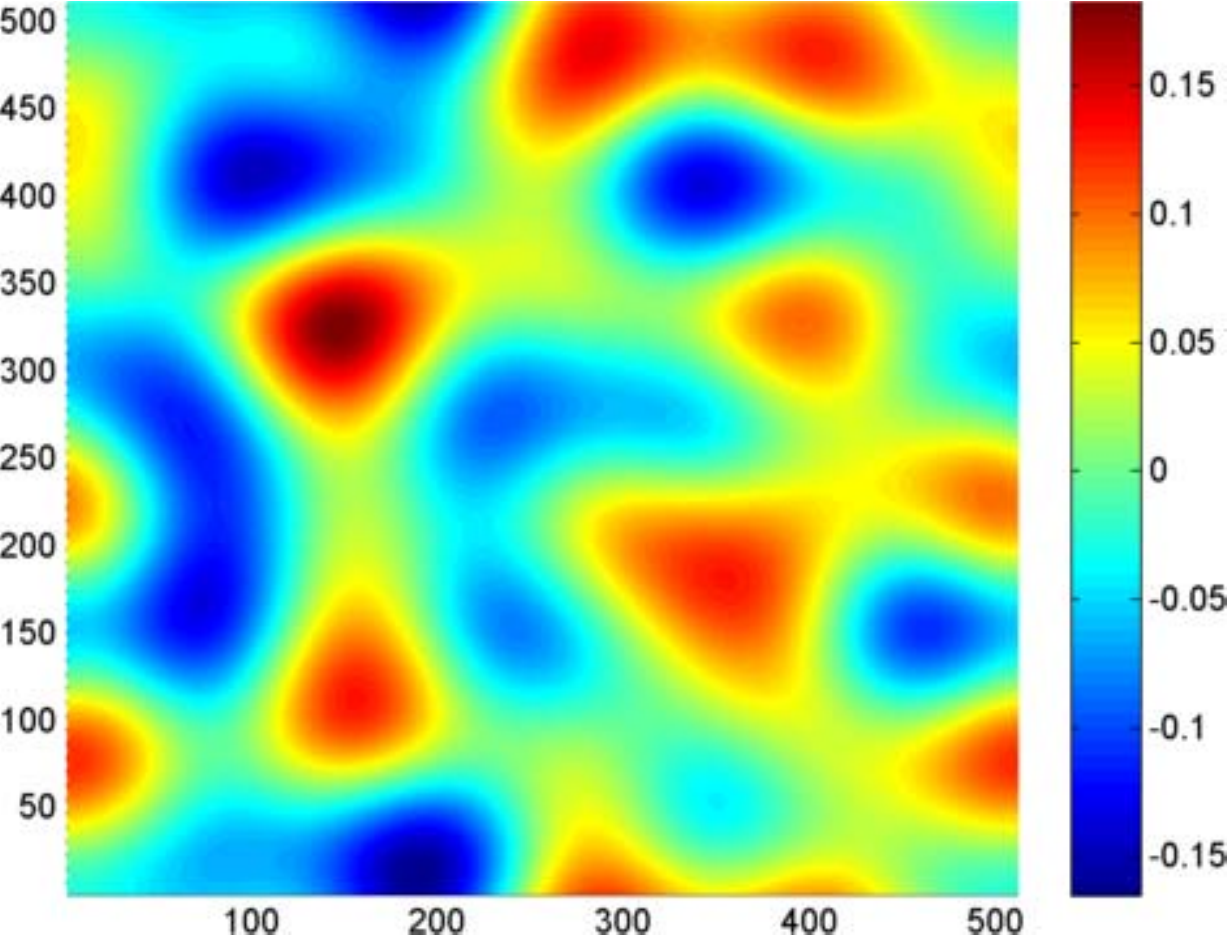}}
\subfloat[$\km=0.1$]{\label{fig:BesLom_2}
\includegraphics[angle=0, width=0.25\linewidth]{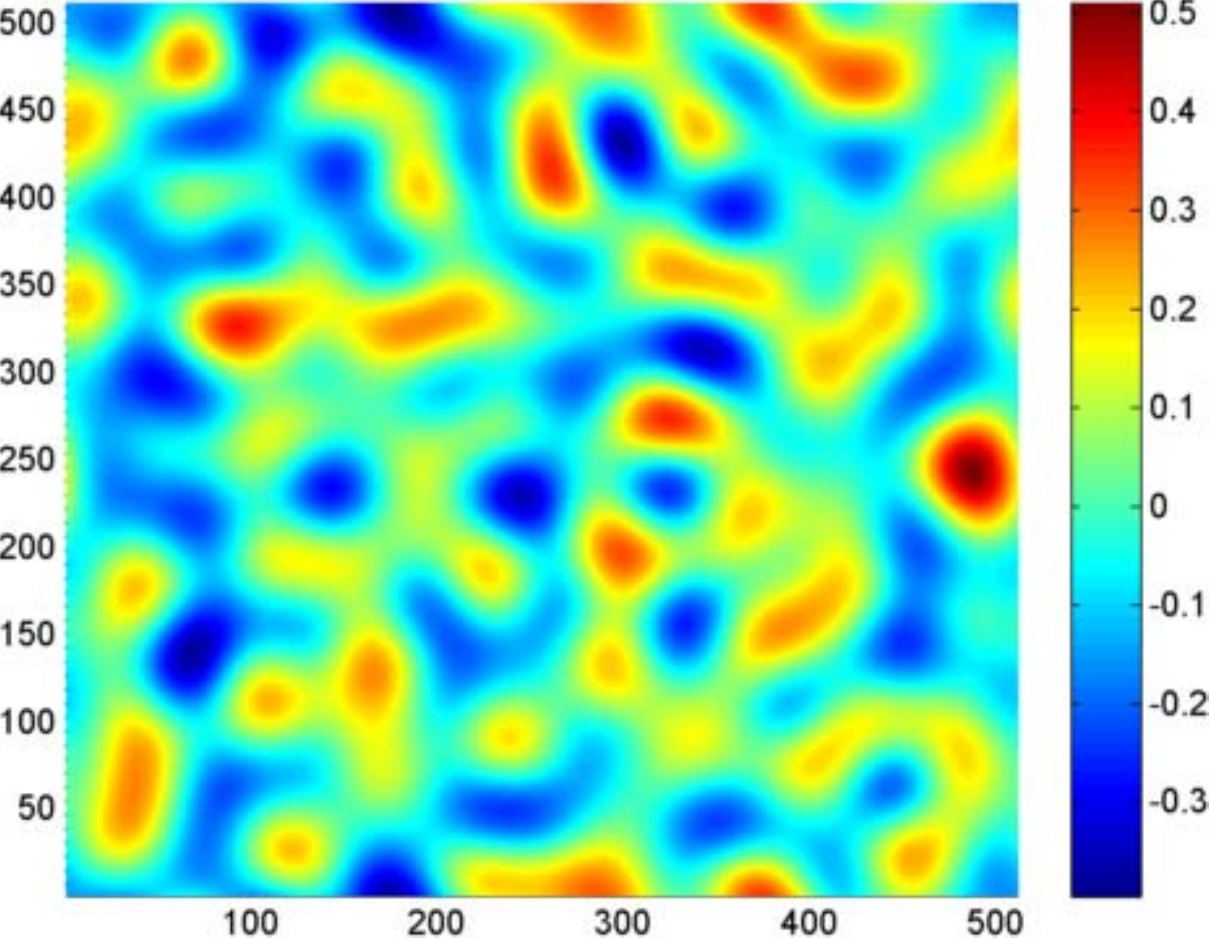}}
\subfloat[$\km=0.2$]{\label{fig:BesLom_3}
\includegraphics[angle=0, width=0.25\linewidth]{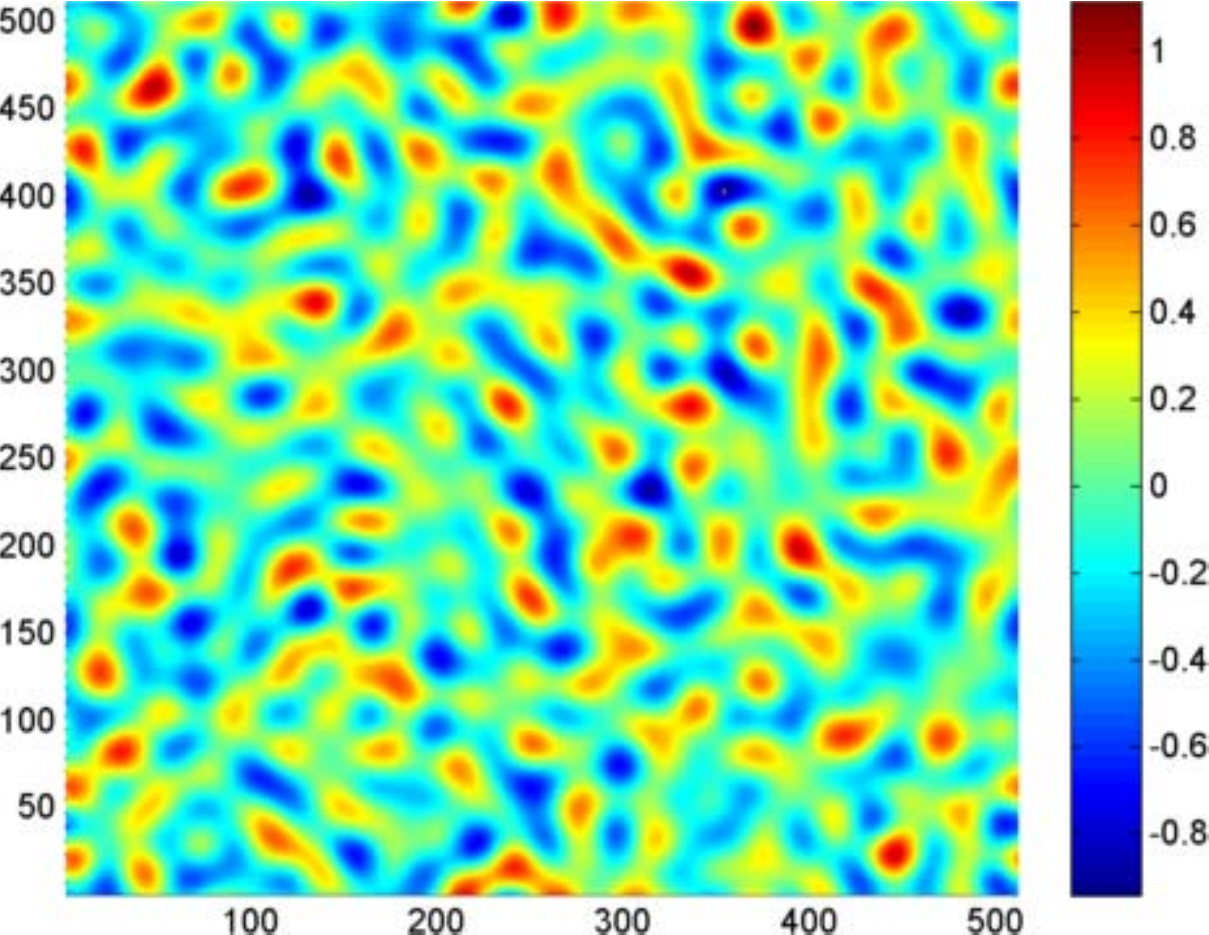}}
\subfloat[$\km=0.5$]{\label{fig:BesLom_4}
\includegraphics[angle=0, width=0.25\linewidth]{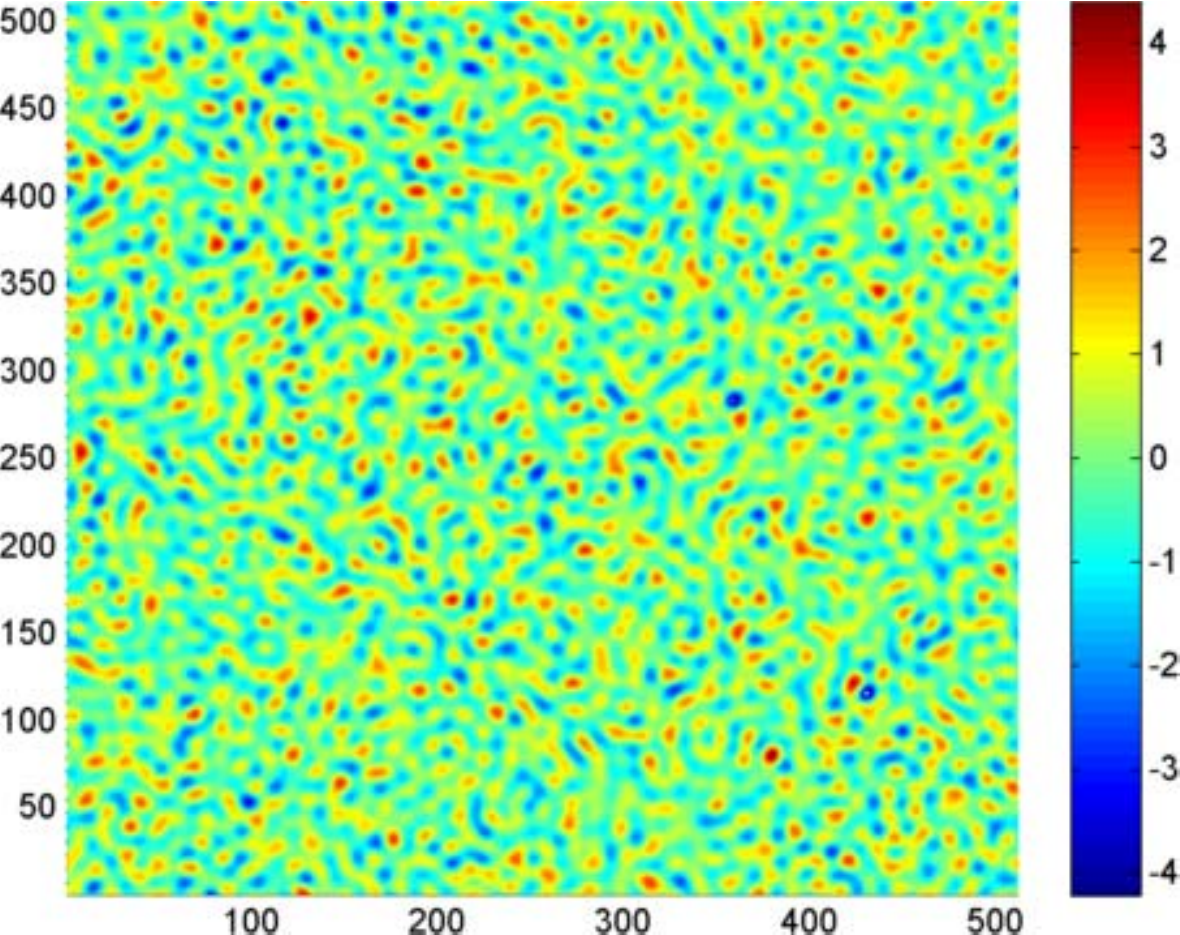}}
\caption{Bessel-Lommel SRF realizations with covariance function~\eqref{eq:cov2d-ssrf} on square $512 \times 512$ grid.  For all realizations  $\eta_{0}=0.01, \eta_{1}=1, \xi=2$.
 FFT spectral simulation with identical random generator seed is used. }
 \label{fig:BesLom-simul}
\end{figure*}

Four different realizations of Bessel-Lommel SRFs are shown
in Fig.~\ref{fig:BesLom-simul}. The realizations have identical $\eta_0$, $\eta_1$ and $\xi$ but different $\km$.
Two trends are obvious in the plots: (i) the variance of the fluctuations increases with $\km$ and (ii)   the
size of characteristic spatial patterns (i.e., areas containing values above or below a given threshold)
is reduced as $\km$ increases. Trend (i) is due to increased spectral weight in the tail of the spectral density function~\eqref{eq:spd-BL}
as $\km \uparrow$.
 Trend (ii) reflects the  decline of the integral range with $\km \uparrow$ according to~\eqref{eq:BL-int-range}
 and as shown in Fig.~\ref{fig:bl-cora-kc}.

\section{Results, Discussion and Conclusions}
\label{sec:concl}
The main results of this paper consist of Propositions 1-4 and Definition~\ref{defi:lac}.
More specifically, we introduce covariance functions with enhanced parametric dependence
motivated by local interaction SRF models.
We derive explicit expressions~\eqref{eq:cov2d} for two-dimensional,
three-parameter, radial SSRF covariance functions at the limit of infinite spectral
cutoff. Gaussian SRFs with such covariance functions are continuous but non-differentiable in the mean square sense.

In addition to SSRF covariance functions, we develop a family of four-parameter, radial Bessel-Lommel covariance
 functions~\eqref{eq:cov-BL} that are valid in $d \ge 2$.
 In contrast with their SSRF counterparts, the Bessel-Lommel functions are infinitely differentiable at the origin, signifying that a Gaussian SRF with such covariance dependence is infinitely differentiable in the mean square sense.

 Since covariance functions are the key ingredient in best linear unbiased estimation methods, e.g., kriging~\cite{Stein99,Rasmussen06},
  the SSRF and Bessel-Lommel functions will provide increased flexibility in the interpolation and simulation
of spatial processes~\cite{Christakos92,Cressie93,Stein99}. In addition, they can be used as  background error covariance models in
  variational assimilation of geophysical data~\cite{Barker04,Weaver13}. They also provide new choices of radial functions suitable for applications in radial basis function interpolation~\cite{Wendland05}.
   Finally, they constitute non-negative definite kernel functions with potential for applications in  machine learning~\cite{Genton02}.

We also demonstrate that the length scales  of SRFs with SSRF and Bessel-Lommel covariance
functions are not uniquely determined by a single characteristic length $\xi$.
We illustrate this behavior by investigating the integral range  of the above covariance models.
In addition, we introduce a correlation spectrum  based on the Fourier transform
of the covariance function's fractional Laplacian of  order   $0 \le \alpha \le 1$. This spectrum allows
quantifying  correlation properties of  mean-square continuous SRFs with a radial, unimodal spectral density.
Finally, we derive explicit expressions for the $\alpha$-spectrum of SSRF and Bessel-Lommel covariance functions.

The assumption of statistical isotropy was used herein primarily for reasons of conciseness
in presentation. Nonetheless, it is straightforward to construct covariance models with geometric (elliptic) anisotropy
 by rotation and rescaling transformations of the coordinate axes~\cite{Gelhar83,dth02,dth08}.


\appendices

\section{Proof of Proposition~\ref{prop:ssrf}}
\label{App:proof-ssrf}
\renewcommand{\theequation}{A-\arabic{equation}}
\renewcommand{\thesubsection}{A-\arabic{subsection}}
\setcounter{subsection}{0}
\setcounter{equation}{0}

\begin{IEEEproof}
Defining  dimensionless wavevectors $u= \kk\,\xi$ and lag distances $h=\rr/\xi$,
 the spectral integral~\eqref{eq:cov2d-spectral} is
simplified as follows:
\begin{equation}
\label{eq:cov2d-spectral-1}
 \Gxx(h;\bmthe) = \frac{\eta_0 }{2\pi} \, \int\limits_{0}^{\infty} du \,\frac{ u \, J_{0}(u\,h)
  }{1+\e u^2+ u^4}.
\end{equation}
Equation~\eqref{eq:cov2d-spectral-1} shows that the only non-trivial parameter is $\e$;
$\eta_0$ is a multiplicative scale factor, whereas the characteristic length $\xi$ is absorbed
in the non-dimensional lag $h.$
The rational function $1/\Pi(u)$, where $\Pi(u)$ is the SSRF characteristic polynomial
defined in~\eqref{eq:char-pol}, admits the following expansion
\begin{equation}
\label{eq:invPI}
\frac{1}{\Pi(u)}= \left\{ \begin{array}{cc}
                    \frac{1}{t_{+}^{\ast}-t_{-}^{\ast}} \, \left( \frac{1}{u^2-t_{+}^{\ast}} -
\frac{1}{u^2 -t_{-}^{\ast}}   \right), & \e \neq 2 \\
                    \frac{1}{(u^2+1)^2} & \e = 2,
                  \end{array} \right.
\end{equation}
where $t_{\pm}^{\ast}={\left(-\e \pm \Delta \right)}/{2}$
are the roots of $\Pi(t=u^2)$.

In light of~\eqref{eq:invPI}, the  integral~\eqref{eq:cov2d-spectral-1} is evaluated
using the Hankel-Nicholson formula~(11.4.44) in~\cite[p. 364]{Abramowitz72}:
\begin{align}
\label{eq:Nicholson}
& \int_{0}^{\infty} du \frac{u^{\nu+1}J_{\nu}(h\,u)}{(u^2 + z^2)^{\mu+1}}  = \frac{h^{\mu}z^{\nu-\mu}}{2^{\mu}\Gamma(\mu+1)}
K_{\nu-\mu}(h\,z).
\end{align}
This equation is valid for  $ h>0, \Re(z)>0$, and $ -1<\Re(\nu)< 2\Re(\mu)+\frac{3}{2}$.
The above is applied to~\eqref{eq:cov2d-spectral-1} with (i) $\e \neq 2$,   $\nu=0$, $\mu=0$, \,
 $z^{2}_{\pm}= - t_{\pm}^{\ast}$  and (ii) $\e = 2$,  $\nu=0$ and $\mu=1$.
 In case (ii) we obtain~\eqref{eq:Cov_eta_eq_2} and in case (i) the following
\begin{equation*}
\Gxx(h ;\bmthe) =\frac{\eta_0 \left[  K_{0} (hz_{+})   -  K_{0} (hz_{-}) \right] }{2\pi \sqrt{\eta_{1}^{2}-4}}
, \; \e \neq 2.
\end{equation*}
The coefficients $z_{\pm}=\sqrt{-t_{\pm}^{\ast}}$
 are plotted versus $\e$ in Fig.~\ref{fig:roots_bess}.
For $\e>2$ both  $z_{+}$ and $z_{-}$
are real numbers, hence proving~\eqref{eq:Cov_eta_gt_2}. For $-2< \e<2$ $\Re(z_{+})=\Re(z_{-})$, whereas
$\Im(z_{+})=-\Im(z_{-})$, i.e., $z_{-} =\overline{z_{+}}$.
The analytic continuation property
$K_{0}(\overline{z})=\overline{K_{0}(z)}$~\cite[p. 377]{Abramowitz72} leads to~\eqref{eq:Cov_eta_lt_2}
 which is explicitly real-valued.

\subsection{Continuity}

A stationary SRF  is mean square continuous $\forall \bfs \in \Rr^d$ if and only
if its covariance function is continuous at zero lag
~\cite{Adler81,Abrahamsen97}. This condition is satisfied for the SSRF covariance.

\subsection{Differentiability}
 Differentiability of the SRF $X(\bfs,\om)$ in the mean-square sense requires that all second-order partial
 derivatives of the covariance function at $\|\bfr \|=0$ exist~\cite[p.~27]{Adler81}.
 This requirement is equivalent to the convergence of the second-order spectral moment
 \[
 \Lambda_{d}^{(2)}:=\int_{\Rr^d} d\bfk \, \kk^2 \, \widetilde{\Gxx}(\bfk,\bmthe).
 \]
For the SSRF spectral density in $d=2$  the above becomes
 \[
 \Lambda_{2}^{(2)} \propto \lim_{\km \ra \infty} \int_{0}^{\km} d\kk \, \frac{\kk^3 }{1 + \e (\kk\xi)^2 + (\kk \xi)^4}.
 \]
This integral develops a logarithmic divergence as $\km \ra \infty$. Hence, the
SSRF is non-differentiable in the mean-square sense.
 \begin{figure}[!t]
  \centering
    \subfloat[Real part of $z_{\pm}=\sqrt{-t_{\pm}^{\ast}}$]{\label{fig:real_r1}
    \includegraphics[width=0.49\linewidth]{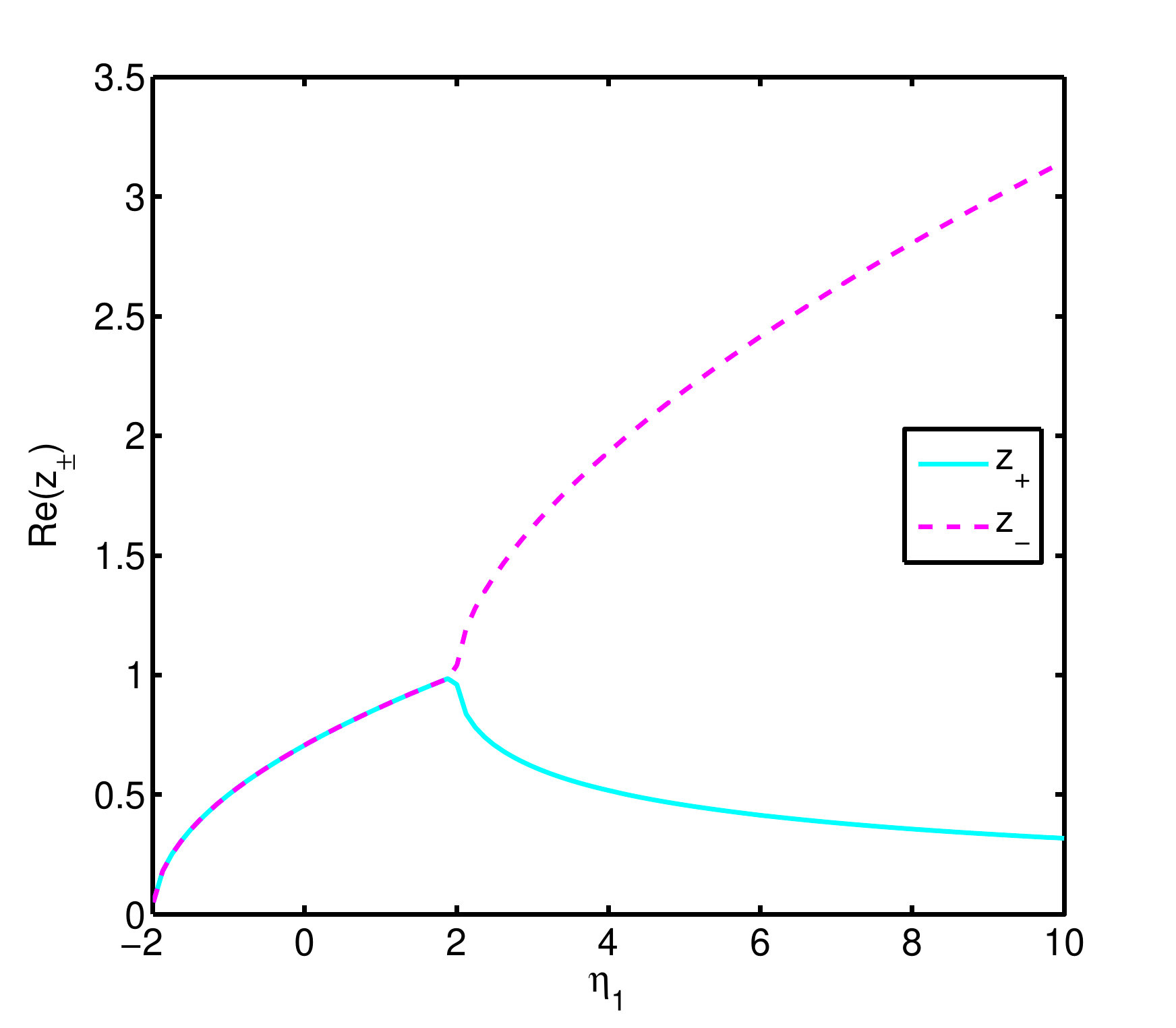}}
    \subfloat[Imaginary part of $z_{\pm}=\sqrt{-t_{\pm}^{\ast}}$]{\label{fig:imag_r1}
    \includegraphics[width=0.49\linewidth]{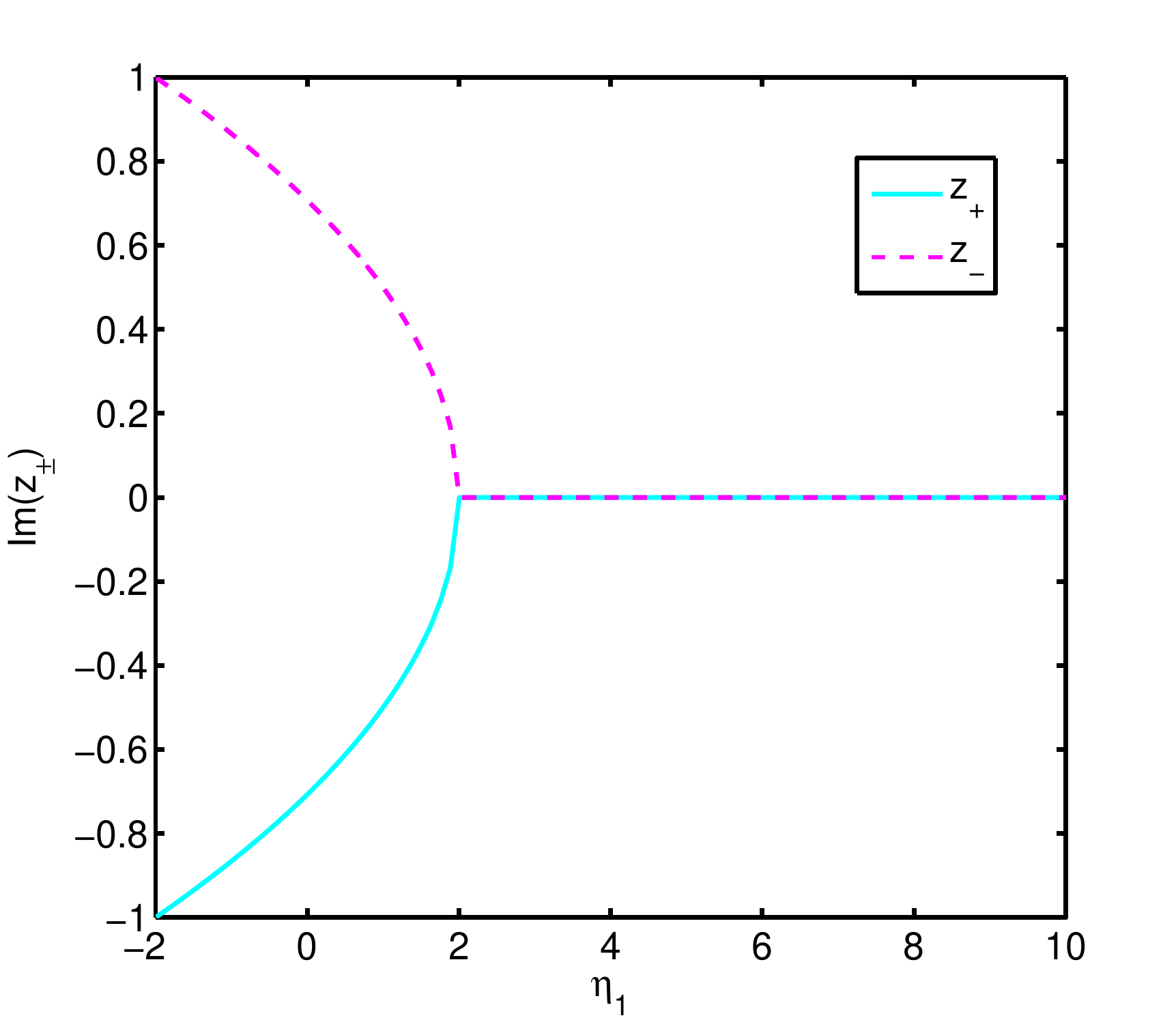}}
  \caption{Real and imaginary parts of the roots $z_{+}=\sqrt{-t_{+}^{\ast}}$ and $z_{-}=\sqrt{-t_{-}^{\ast}}$
  of the characteristic polynomial $\Pi(u)=1+\eta_1 u^2 + u^4$ according to~\eqref{eq:zpm}.}
  \label{fig:roots_bess}
\end{figure}

 \hspace{\stretch{1}}
\end{IEEEproof}

\section{Proof of Proposition~\ref{prop:BL}}
\label{App:proof-BL}
\renewcommand{\theequation}{B-\arabic{equation}}
\renewcommand{\thesubsection}{B-\arabic{subsection}}
\setcounter{subsection}{0}
\setcounter{equation}{0}

\begin{IEEEproof}
Let $z=u_{c}\, h$,  $\nu=d/2-1$,  $\mu > -(\nu+1)$, and $J_{\nu}(\cdot)$ be the
Bessel function of the first kind of order $\nu$. We then define the function
\begin{equation}
\label{eq:Amn-def}
\A_{\mu,\nu}(z)=\int_{0}^{1} dx \, x^{\mu} J_{\nu}(z\,x).
\end{equation}
Then, $\A_{\mu,\nu}(z)$ is evaluated using~\cite[p.~676, eq.~(6.561.13)]{Gradshteyn07} as follows
\begin{align}
\label{eq:Amn}
\A_{\mu,\nu}(z)  = &
\frac{2^{\mu} \Gamma\left( \frac{\nu+\mu+1}{2}\right)}{z^{\mu+1} \, \Gamma\left( \frac{\nu-\mu+1}{2}\right)}
+ \frac{(\mu + \nu -1) J_{\nu}(z) S_{\mu-1,\nu-1}(z)}{z^{\mu}}
\nonumber \\
&   \quad \quad \quad  -  \frac{J_{\nu-1}(z) S_{\mu,\nu}(z)}{z^{\mu}}.
\end{align}

 We use the normalizing variable transformations $x=\kk/\km$, $h=\rr/\xi,$ and $u_{c}=\km \xi$.
In view of the dimensionless variables $x, h, u_c$, the integral~\eqref{eq:cov-BL-spectral} becomes
\begin{align}
\label{eq:cov-BL-spectral-1}
\inc^{BL}(h;\bmthe)  =  &   \frac{ u_c^{1+d/2}\, h^{1-d/2}}{(2\pi)^{d/2} \eta _0\, \xi^{2d} } \int_{0}^{1} d x \,
 {x^{d / 2} J_{d/2-1}(x h u_c)}
 \nonumber \\
 &  \quad  \cdot \left[ 1+\e (x u_c)^2+ (x u_c)^4 \right].
\end{align}

In light of~\eqref{eq:Amn-def}, \eqref{eq:cov-BL-spectral-1} and using
$z =u_{c} h = \km \rr$ as the dimensionless distance, the function $\inc^{BL}(\bfr;\bmthe) $
defined by~\eqref{eq:cov-BL-spectral} is given by
\begin{align}
\label{eq:cov-BL-spectral-2}
\inc^{BL}(z;\bmthe) = &  \frac{g_{0}(\bmthe)}{ z^{\nu}} \left[\A_{\nu+1,\nu}(z)
+ \e u_{c}^2 \A_{\nu+3,\nu}(z) \right.
\nonumber \\
    &  \left. \quad \quad \quad + u_{c}^4 \A_{\nu+5,\nu}(z)\right] ,
\end{align}
\begin{equation}
\label{eq:c}
g_{0}(\bmthe) = \frac{ \km^{d}}{(2\pi)^{d/2} \eta _0\, \xi^{d} }.
\end{equation}

For the three terms $\A_{\mu,\nu}(z)$ $(\mu=\nu+1, \nu+3, \nu +5)$ included in $\inc^{BL}(z;\bmthe)$, the
parameters $\mu, \nu$ satisfy the relation
\begin{equation}
\label{eq:nu-m-mu}
\nu - \mu +1
= - 2\,l, \;\mbox{where} \;  l=0,1,2.
\end{equation}
Equations~\eqref{eq:cov-BL} follow directly from~\eqref{eq:cov-BL-spectral-2} which expresses $\inc^{BL}(z;\bmthe)$
in terms of $\A_{\mu,\nu}(z)$, and from~\eqref{eq:Amn} which expresses the integrals $\A_{\mu,\nu}(z)$ in terms of
Lommel functions.  In view of~\eqref{eq:nu-m-mu}, the Gamma function
contributions to  $\A_{\nu+2l+1,\nu}(z)$  in~\eqref{eq:Amn}
vanish due to the poles of $\Gamma(n)$ at $n \in \mathbb{Z}_{0,-}$.

\subsection{Permissibility}
The non-negative definiteness of $\inc^{BL}(z;\bmthe)$ is based on
Bochner's theorem and
the fact that, according to~\eqref{eq:spd-BL}, $\widetilde{\inc^{BL}}(\kk;\bmthe) \ge 0$  for $\e>-2$.

\subsection{Differentiability}
The existence of the $n$-th order partial derivatives of the Bessel-Lommel SRF in the mean-square sense requires that
all the partial derivatives of order $2n$ of $\inc^{BL}(z;\bmthe)$ exist at $z=0$. This condition  is ensured by the convergence of the $2n$-th order spectral moment
 \[
\Lambda_{d}^{(2n)}  = \int_{0}^{\km } d\kk \, {\kk^{n+d-1} }{\, \left[ 1 + \e (\kk\xi)^2 + (\kk \xi)^4 \right]}.
 \]

\hspace{\stretch{1}}
\end{IEEEproof}
%

\section{Proof of Proposition~\ref{prop:ssrf-scales}}
\label{App:Proof-C}
\renewcommand{\theequation}{C-\arabic{equation}}
\renewcommand{\thesubsection}{C-\arabic{subsection}}
\setcounter{subsection}{0}
\setcounter{equation}{0}

\begin{IEEEproof}
To find the supremum of $f(\kk):=\kk^{2\alpha} \, \widetilde{\Gxx}(\kk;\bmthe)$ we consider the extremum
condition ${\rm d} f(\kk)/{\rm d} \kk =0$, which admits the following two roots:
\[
\tilde{\kappa}_{1,2} = \sqrt{\frac{ \pm \sqrt{\e^2 \, (1 - \alpha)^2 - 4 \alpha (\alpha-2)} - \e \, (1 - \alpha)}{2(2 - \alpha) \xi^2}}.
\]
For $ 0 \le \alpha < 1$ only $\tilde{\kappa}_{1} \in \Rr$ and  $\sup f(\kk) = f(\tilde{\kappa}_{1})$.

According to~\eqref{eq:ssrf-spd} the denominator of~\eqref{eq:micro} becomes
 \begin{align}
\label{eq:micro-ssrf-den}  &
\mathcal{S}_d \, \eta_{0} \, \xi^{1-2\alpha} \, \int_{0}^{\infty} dx \,\phi_{\alpha}(x)
\\
\label{eq:phia}
 \mbox{where} \quad \phi_{\alpha}(x)   = &     \frac{x^{1+2\alpha}}{ 1 + \e\, \xi^2 \, x^2 + \xi^4 \, x^4 }.
\end{align}
To simplify notation we define
\begin{equation}
\label{eq:def-Ia}
I_{\alpha}(\phi):=\int_{0}^{\infty} dx \phi_{\alpha}(x).
\end{equation}
In order to calculate the integral~\eqref{eq:micro-ssrf-den} we use \emph{Lebesgue's  dominated convergence
theorem}~\cite{Schwartz08} expressed as follows:

\begin{theorem}
\label{thm:dom-conv}
Let $\phi_{\alpha}(x)$ be a real-valued function $\forall x \in \Rr$ which is integrable $\forall \alpha \in [0, 1]$. If there is a real-valued function $g_{n}(x)$ such that (i)
$\lim_{n \ra \infty} \phi_{\alpha}(x) \, g_{n}(x) = \phi_{\alpha}(x), \, \forall x \in \Rr$  and (ii) $ |\phi_{\alpha}(x) \, g_{n}(x)| \le \phi^{\ast}(x), \forall x \in \Rr$,
where $\phi^{\ast}(x)$ is an integrable function, then
\begin{align*}
 I_{\alpha}(\phi)  & = \int_{0}^{\infty} dx \lim_{n \ra \infty} g_{n}(x)  \phi_{\alpha}(x)
    = \lim_{n \ra \infty} \int_{0}^{\infty} dx g_{n}(x)  \phi_{\alpha}(x).
\end{align*}
\end{theorem}

We define the following auxiliary function
\begin{equation}
\label{eq:gn}
g_{n}(x) : = \frac{ 2^{\alpha }\, \Gamma(\alpha+1) \, J_{\alpha}(x/n)}{ \left(\frac{x}{n}\right)^\alpha }.
\end{equation}
Condition (i) of Theorem~\ref{thm:dom-conv} is satisfied because $\lim_{n \ra \infty}  g_{n}(x)=1$
based on the infinite series expansion of the Bessel function of the first kind around zero~\cite[p.~40]{Watson95}.

To prove the condition (ii) we apply the following steps.
\begin{enumerate}
\item
 For condition (ii) it suffices that $ |\phi_{\alpha}(x) \, g_{n}(x)| \le \phi_{\alpha}(x) $, because $\phi_{\alpha}(x)$ given by~\eqref{eq:phia}
is integrable.
Given that  $\phi_{\alpha}(x) >0$, it suffices to show that
$| g_{n}(x) |  \le 1$.
\item
We use the integral representation of $J_{\alpha}(z)$ given by~\cite[8.411.4]{Gradshteyn07}, where $z \in \Rr$:
\[
J_{\alpha}(z) = \frac{2\left( \frac{z}{2} \right)^{\alpha}}{\Gamma(\alpha+1/2)\Gamma(1/2)} \,
\int_{0}^{\pi/2} d\theta \sin^{2\alpha}\theta \cos \left( z \cos\theta \right)
\]
\item
Since $|\sin^{2\alpha}(\theta) \cos \left( z \cos\theta\right)| \le 1$ and $\Gamma(1/2) = \sqrt{\pi}$
it follows from the above
that $|J_{\alpha}(z)| \le \frac{\left( \frac{z}{2} \right)^{\alpha} \sqrt{\pi}}{\Gamma(\alpha+1/2)}.$
\item
In light of this inequality and~\eqref{eq:gn},  proving that $|g_{n}(x)| \le 1$ is equivalent to showing that
$ \mu_{\alpha}:= \Gamma(\alpha +1) / \Gamma(\alpha +1/2) \le \sqrt{\pi}.$
\item
Based on the inequality $ \mu_{\alpha} < \sqrt{\alpha +1/2}$ (valid for $\alpha > -1/4$)
 the maximum upper bound of $ \mu_{\alpha}$ for $0 \le \alpha \le 1$ is $\sqrt{3/2} < \sqrt{\pi}$.
 Hence, in light of the previous step $|g_{n}(x)| \le 1$. This concludes the proof of condition (ii).
 \end{enumerate}

In light of the above, we can use dominated convergence to calculate $I_{\alpha}(\phi)$ as follows
\begin{equation}
\label{eq:Ia-dc}
I_{\alpha}(\phi) = \lim_{n \ra \infty}  (2 n)^{\alpha } \Gamma(\alpha+1)   \,\tilde{I}_{\alpha}(\phi)
\end{equation}
where
\begin{equation}
\label{eq:tIa-dc}
\tilde{I}_{\alpha}(\phi) = \int_{0}^{\infty} dx \,
\frac{J_{\alpha}(x/n)   x^{1+\alpha}}{ 1 + \e \xi^2  x^2 + \xi^4  x^4 }.
\end{equation}
The  integral $\tilde{I}_{\alpha}(\phi)$ is evaluated by means of the Hankel-Nicholson formula~\eqref{eq:Nicholson} ($\nu =\alpha,$ $\mu=0$ for $\e \neq 2$ , $\mu=1$ for $\e= 2$) which leads to
\begin{align}
\label{eq:tIa-dc-1}
\tilde{I}_{\alpha}(\phi) =  \left\{
                                    \begin{array}{cc}
                                     \frac{z_{+}^{\alpha} \, K_{\alpha}(z_{+}/n) - z_{-}^{\alpha} \,K_{\alpha}(z_{-}/n)}{\sqrt{\e^2-4}}  &    \e \neq 2\\
                                     \frac{K_{\alpha-1}(1/n)}{2n} = \frac{K_{1-\alpha}(1/n)}{2n}       &    \e = 2.
                                   \end{array}
                                   \right.
\end{align}

To evaluate   $\lim_{n \ra \infty}\tilde{I}_{\alpha}(\phi)$ for $1> p >0$
we use the series expansion~\cite{Schwartz08} of the K-Bessel function
\begin{align*}
K_{p}(x) = &    \frac{1}{2} \left[ \Gamma(p) \left( \frac{2}{x}\right)^{p} \left( 1+ O(x^2) \right) +
 \Gamma(-p) \left( \frac{x}{2}\right)^{p} \right.
 \nonumber \\
 &   \left. \quad \quad \left( 1+ O(x^2) \right) \right].
\end{align*}

For $\e =2$, $p = 1 - \alpha$, and $x=1/n$ the dominant contribution at  $n \ra \infty$ comes from
the $O(x^{-p})$ term of the first series on the right hand side, which gives
$ K_{1-\alpha}(1/n) \approx \frac{1}{2}  \Gamma(1 - \alpha) (2n)^{1 - \alpha}$.

For $\e \neq 2$  the $O(x^{-p})$ term of the first series
cancels out due to the difference between the two Bessel functions, whereas the $O(x^{2-p})$ terms
vanish at the limit $n \ra \infty$.
A finite contribution
comes from the  $O(x^p)$ term of the second series on the right hand side, i.e.,
 $z_{+}^{\alpha} \, K_{\alpha}(z_{+}/n) - z_{-}^{\alpha} \,K_{\alpha}(z_{-}/n) \sim
 \frac{\Gamma(-\alpha)}{2(2n)^{\alpha}}\left( z_{+}^{2\alpha} -  z_{-}^{2\alpha} \right) $.

Thus, based on the above asymptotic analysis of the K-Bessel function,~\eqref{eq:Ia-dc},~\eqref{eq:tIa-dc}, and~\eqref{eq:tIa-dc-1} we
obtain the following equation (where $\Delta = \sqrt{\e^2-4}$):
\begin{align}
\label{eq:Ia}
I_{\alpha}(\phi) = &  \left\{
\begin{array}{cc}
\frac{\Gamma(1-\alpha)\Gamma(1+\alpha) \, \left[\left( \e + \Delta \right)^{\alpha}  -
\left( \e - \Delta\right)^{\alpha}\right]}{2^{\alpha+1}\alpha \Delta}
   &    \e \neq 2\\
   \frac{\Gamma(1-\alpha) \, \Gamma(1+\alpha) }{2}       &    \e = 2.
   \end{array}
   \right.
\end{align}

Finally, based on~\eqref{eq:Ia}, the definition~\eqref{eq:def-Ia} and~\eqref{eq:micro-ssrf-den},~\eqref{eq:micro-ssrf} is proved.

 \hspace{\stretch{1}}
\end{IEEEproof}

\section{Proof of Proposition~\ref{prop:BL-scales}}
\label{App:Proof-D}
\renewcommand{\theequation}{D-\arabic{equation}}
\renewcommand{\thesubsection}{D-\arabic{subsection}}
\setcounter{subsection}{0}
\setcounter{equation}{0}

\begin{IEEEproof}

Based on the spectral density~\eqref{eq:spd-BL} it follows that the denominator in~\eqref{eq:micro} is given by
\begin{align}
\label{eq:BL-micro-den}
\int_{\Rr^d} d\bfk \,\kk^{2\alpha} & \widetilde{\Gxx}^{BL}(\kk;\bmthe) =  \frac{\mathcal{S}_{d}  \km^{d+2\alpha}}{\eta_{0}  \xi^d}
\left( \frac{1}{d+2\alpha} \right.
\nonumber \\
    &   \left. + \frac{\e \km^2 \xi^2}{d+2\alpha+2}  + \frac{ \km^4 \xi^4}{d+2\alpha+4} \right).
\end{align}

Let us define the function $\phi(\kk) = \kk^{2\alpha} \,\widetilde{\Gxx}^{BL}(\kk;\bmthe)$.
The numerator in~\eqref{eq:micro} is then given by
$\sup_{\kk \in \Rr} \, \phi(\kk) = \phi(\kappa^{\ast})$ where $\kappa^{\ast} = \underset{{\kk} \in [0, \km]}{\arg \max} \phi(\kk)$.

\subsection{Non-negative $\e$}
For $\e  \ge 0$, $\phi(\kk)$ is a monotonically increasing function of $\kk$; thus, $\kappa^{\ast}= \km$ and
$ \phi(\kappa^{\ast}) =  \frac{\km^{2\alpha}}{\eta_{0} \xi^d}\, \left( 1 + \e\, \km^2 \, \xi^2 + \km^4 \xi^4  \right)$.
In light of~\eqref{eq:BL-micro-den}, this leads to~\eqref{eq:BL-micro-pos-eta1}.

\subsection{Negative $\e$}
 For $\e <0$, $\phi(\kk)$ develops local extrema at the wavenumbers that solve the equation $d\phi(\kk)/d\kk =0$, i.e., at the~$\kappa_{\pm}$ given by~\eqref{eq:BL-kappa}~\footnote{There are two additional solutions of opposite sign than $\kappa_{\pm}$ which are not further considered, since they are either complex or negative real numbers.}.

\paragraph{Complex $\kappa_{\pm}$}
For $\e^2 < 4 \alpha (\alpha +2) / (\alpha +1)^2$, the $\kappa_{\pm}$ are complex numbers and thus
$\phi(\kk)$ does not have local extrema for $\kk \in \Rr$. Hence, $\kappa^{\ast} = \km$ and
$\lambda_{c}^{(\alpha)}$ is given by~\eqref{eq:BL-micro-pos-eta1}.

\paragraph{Real $\kappa_{\pm}$}
For  $4 > \e^2 \ge 4 \alpha (\alpha +2) / (\alpha +1)^2$,
the $\kappa_{\pm}$ are real numbers,
one corresponding to the position of a local minimum whereas the other to a local maximum.

 If $\alpha>0$, $d\phi(\kk)/d\kk \propto 2\alpha\, (\kk \xi)^{2\alpha-1}$
for  $\kk \ll 1$,   and thus $\phi(\kk) \uparrow$.   Hence,
the maximum of $\phi(\kk)$ occurs at lower $\kk$  than the minimum, i.e., at $\kappa_{-} < \kappa_{+}$.
If $\alpha=0$, then $d\phi(\kk)/d\kk \propto -2|\e|\, \kk \,\xi$ and thus $\phi(\kk) \downarrow$
 for  $\kk \ll 1$.   In this case as well
$\phi(\kk)$ reaches the maximum at  $\kappa_{-}=0,$ whereas the minimum occurs at
$\kappa_{+} = \sqrt{|\e|/2} / (\km \xi)$.
Again, we distinguish between two cases:

\begin{enumerate}

\item If $\kappa_{-} > \km$, then  $\kappa^{\ast} = \km$ and
$\lambda_{c}^{(\alpha)}$ is given by~\eqref{eq:BL-micro-pos-eta1}.
 For $\alpha=0$ it holds that $\kappa_{-}=0$ and thus $\kappa_{-} < \km$.

\item $\km > \kappa_{-}$.
We further distinguish the following cases.

\begin{enumerate}
\item If $\km < \kappa_{+}$, then
$\kappa^\ast = \kappa_{-}$ and $\lambda_{c}^{(\alpha)}$ is given by~\eqref{eq:BL-micro-neg-eta1-1}.

\item If $\km > \kappa_{+}$, then
$ \kappa^{\ast} =  {\arg\max} (\phi(\kappa_{-}), \phi(\km))$ and $\lambda_{c}^{(\alpha)}$ is given by~\eqref{eq:BL-micro-neg-eta1-2}.
\end{enumerate}

\end{enumerate}

\hspace{\stretch{1}}
\end{IEEEproof}


\section*{Acknowledgment}
The author acknowledges funding from the project SPARTA 1591: ``Development of
Space-Time Random Fields based on Local Interaction Models and Applications in the Processing of
Spatiotemporal Datasets''. SPARTA is implemented under the ``ARISTEIA'' Action of the
 operational programme ``Education and Lifelong Learning'' and is co-funded by the European Social Fund
 and National Resources.

\bibliography{IEEEabrv}


\onecolumn
\section{Notation}
\label{app:table}

\begin{small}
\begin{center}
\begin{longtable}{lll}
   \hline
  Symbol & Definition \\
      \hline    \hline
\multicolumn{2}{c}{{\em Direct and Reciprocal Space}} \\
\hline
   $d$ & Spatial dimension\\
   $\D$ & Spatial domain $ \D \in \Rr^d$ \\
   $\bfs \in \Rr^{d}$ & Position vector \\
   $\bfr \in \Rr^{d}$ & Lag vector \\
   $\rr \in \Rr$  & Euclidean norm of $\bfr$ \\
   $\bfk \in \Rr^{d}$ & Wavevector in reciprocal space \\
   $\kk \in \Rr$  & Euclidean norm of $\bfk$ \\
   $ h$ & Dimensionless lag: $h=\rr/\xi$ \\
   $ z$ & Dimensionless lag: $z=\km \, \rr$ \\
   $ \mathcal{S}_d $ &  Surface of unit sphere in $\Rr^{d}$: $ \mathcal{S}_d = 2 \pi^{d/2}/\Gamma(d/2)$ \\
\hline
\multicolumn{2}{c}{{\em Random Fields}} \\
\hline
    $\Omega$  & Event space \\
    $\om$  & State index \\
   $X(\bfs,\om)$ & Scalar spatial random field (SRF) \\
   $x(\bfs)$ & SRF sample function (realization) \\
   $\Gxx(\bfr;\bmthe)$ & Covariance function \\
   $\sigma^{2}_{\rmx}$ & Variance \\
   $\ell_c$ & Integral range of isotropic covariance\\
   $r_c$ & Correlation radius of isotropic covariance\\
   $\alpha$ & Fractional exponent $0 \le \alpha \le 1$ \\
   $\lambda_c^{(\alpha)}$ & Spectrum of SRF length scales\\
\hline
\multicolumn{2}{c}{{\em Spectral Transforms}} \\
\hline
\vspace{0.5pt}\\
    $\widetilde{\Gxx}(\bfk;\bmthe)$ & Covariance spectral density \\
   $\mathbbm{1}_{A}(x)$ & Indicator function of set $A$\\
   $\Fr[ \cdot]$ &  Fourier transform \\
   $\iFr[ \cdot]$ &  Inverse Fourier transform \\
   $\delta(\bfs-\bfs')$ &  Delta function \\
\hline
\multicolumn{2}{c}{{\em FGC SSRF}} \\
\hline
   $\eta_0$ & SSRF amplitude coefficient \\
   $\e$ & SSRF rigidity coefficient \\
   $\xi$ & SSRF characteristic length \\
   $\km$ & SSRF wavevector cutoff \\
   $u_c$ & Dimensionless wavevector cutoff: $u_c =\km \xi$  \\
   $\bmthe$ & Vector of SSRF parameters: $\bmthe=(\eta_0,\e,\xi,\km)^T$ \\
   $\bmthe'$ & Reduced vector of SSRF parameters: $\bmthe'=(\e,\xi,\km)^T$  \\
   $\Pi(\kk \xi)$ & FGC-SSRF characteristic polynomial: $\Pi(u)= 1 + \e\, u^2 +  u^4$  \\
   $t^{\ast}_{\pm}$ & Roots of characteristic polynomial \\
   $ \Delta$ & Discriminant of characteristic polynomial: $\Delta=\sqrt{\e^2-4}$ \\
   \hline
   \multicolumn{2}{c}{{\em Bessel-Lommel SSRF}} \\
\hline
   $ c_{0}$ & Zero-degree polynomial coefficient in Bessel-Lommel spectral density\\
   &     $c_{0}= 1/(\eta_{0} \xi^d)$ \\
   $ c_{1}$ &  Second-degree polynomial coefficient  in Bessel-Lommel spectral density\\
   &     $c_{1}= \e\, \xi^2/(\eta_{0} \xi^d)$ \\
   $ c_{2}$ &  Fourth-degree polynomial coefficient in Bessel-Lommel spectral density\\
   &     $c_{2}=  \xi^4/(\eta_{0} \xi^d)$ \\
   $\widetilde{\inc^{BL}}({\kk};{\bm \theta})$  & Bessel-Lommel spectral density \\
   $\inc^{BL}({\bfs}-{\bfs}';{\bm \theta})$ & Bessel-Lommel covariance kernel \\
   $ g_{0}(\bmthe)$ & Coefficient in Bessel-Lommel  covariance: $g_{0}(\bmthe) = { \km^{d}} /{(2\pi)^{d/2} \eta _0\, \xi^{d} }$ \\
   $ g_{1}(\bmthe)$  & Coefficient in Bessel-Lommel  covariance: $g_{1}(\bmthe)=\e \, u_c^{2} \, g_{0}(\bmthe)$ \\
   $ g_{2}(\bmthe)$  & Coefficient in Bessel-Lommel  covariance: $g_{2}(\bmthe)=u_c^{4} \, g_{0}(\bmthe)$ \\
\hline
\multicolumn{2}{c}{{\em Certain Functions and Operators}} \\
\hline
    $ \Gamma(\nu)$ & Gamma function \\
   $J_{\nu}(x)$ &   Bessel function of the first kind of order $\nu$ \\
   $K_{\nu}(x)$ &   Modified Bessel function of the second kind of order $\nu$ \\
   $S_{\mu,\nu}(z)$ &   Lommel functions \\
   $\A_{\mu,\nu}(z)$ & $\int_{0}^{1} dx \, x^{\mu} J_{\nu}(z\,x)$ \\
   $ \nabla$ &  Gradient operator: $ \nabla = \left( \frac{\partial}{\partial s_{1}}, \ldots,
   \frac{\partial}{\partial s_{d}} \right)^T$ \\
   $ \nabla^2 = \triangle $ &  Laplacian operator: $ \nabla^2 = \sum_{i=1}^{d} \frac{\partial^2}{\partial s_{1}^{2}} + \ldots
   \frac{\partial}{\partial s_{d}^{2}} $ \\
   \hline \hline
\caption{Notation.}
\label{tab:Notation}
\end{longtable}
\end{center}
\end{small}

\twocolumn


\end{document}